\documentclass[final,3p]{elsarticle}
 \usepackage[latin1]{inputenc}
 \usepackage{graphics}
 \usepackage{graphicx}
 \usepackage{epsfig}
\usepackage{amssymb}
 \usepackage{amsthm}
 \usepackage{lineno}
 \usepackage{amsmath}
   \numberwithin{equation}{section}
\usepackage{mathrsfs}

\newtheorem{thm}{Theorem}[section]

\newtheorem{lem}[thm]{Lemma}

\newtheorem{defn}[thm]{Definition}

\biboptions{sort&compress,square}
\begin{document}
\begin{frontmatter}
\author[rvt1]{Jian Wang}
\ead{wangj@tute.edu.cn}
\author[rvt2]{Yong Wang\corref{cor2}}
\ead{wangy581@nenu.edu.cn}

\cortext[cor2]{Corresponding author.}
\address[rvt1]{School of Science, Tianjin University of Technology and Education, Tianjin, 300222, P.R.China}
\address[rvt2]{School of Mathematics and Statistics, Northeast Normal University,
Changchun, 130024, P.R.China}

\title{ Equivariant Twisted Bismut Laplacian with Torsion \\ and  KKW type theorems }
\begin{abstract}
This paper aims to  provide an explicit computation of the noncommutative residue density
associated with equivariant twisted Bismut Laplacian with torsion on compact manifolds with (or without) boundary.
We prove the equivariant twisted Kastler-Kalau-Walze type theorems with torsion
on compact manifolds   with boundary.
\end{abstract}
\end{frontmatter}
\section{Introduction}
\label{1}

Let $E$ be a finite-dimensional complex vector bundle over a closed compact manifold $M$
of dimension $n$, the noncommutative residue of a pseudo-differential operator
$P\in\Psi DO(E)$ \cite{Wo,Wo1,Gu} can be defined by
 \begin{equation}
res(P):=(2\pi)^{-n}\int_{S^{*}M}\mathrm{Tr}(\sigma_{-n}^{P}(x,\xi))\mathrm{d}x \mathrm{d}\xi,
\end{equation}
where $S^{*}M\subset T^{*}M$ denotes the co-sphere bundle on $M$ and
$\sigma_{-n}^{P}$ is the component of order $-n$ of the complete symbol
$
\sigma^{P}:=\sum_{i}\sigma_{i}^{P}
$
of $P$, and the linear functional $res: \Psi DO(E)\rightarrow \mathbb{C }$
is in fact the unique trace (up to multiplication
by constants) on the algebra of pseudo-differential operators $\Psi DO(E)$.
Connes proved that the noncommutative residue on a compact manifold $M$ coincided with Dixmier's trace on pseudodifferential
operators of order -dim$M$, and derived a conformal 4-dimensional Polyakov action analogy by using the noncommutative residue \cite{Co1}.
More precisely, Connes \cite{Co2,Co3} made a challenging observation that
the  noncommutative residue ${\rm  Wres}(D^{-2})$ was proportional to the Einstein-Hilbert action of general relativity
 (which, in turn was essentially the $a_{2}$ coefficient of the heat expansion of Laplacian operator $\Delta$).
This consequence's  earlier proofs of Kastler, Kalau and Walze ( called the Kastler-Kalau-Walze theorem) were based on the explicit computation of
$\sigma_{-4}{(D^{-2})}$, which was most easily performed in normal coordinates (this was the main advantage of the KW approach \cite{KW}
over the ¡®brute force¡¯ computation by Kastler \cite{Ka}).
On the other hand, it was Boutet de Monvel \cite{BDM} who brought in the operator-algebraic aspect
with his calculus established, and constructed a relatively
small `algebra', containing the elliptic differential boundary value problems as well
as their parametrices. In \cite{ES2}, Schrohe provided an introduction to Boutet de Monvel¡¯s calculus on the half-space $\mathbb{R}^{n}_{+}$
in the framework of a pseudodifferential calculus with operator-valued symbols.
For the operators of order $-n$ and type zero in Boutet de Monvel¡¯s calculus,
Nest and Schrohe \cite{NS} showed that Dixmier¡¯s trace can be computed in terms of the same
expressions that determined the noncommutative residue, and provided how to compute Dixmier's trace
for parametrices or inverses of classical elliptic boundary value problems.
Furthermore, Fedosov et al. defined the noncommutative residue on Boutet de Monvel's algebra and proved that it was a
unique continuous trace for the case of manifolds with boundary \cite{FGLS}.
Then Wang generalized Connes' results to the case of manifolds with boundary,
and proved  Kastler-Kalau-Walze type theorems for lower-dimensional manifolds
with boundary associated with Dirac operator and signature operator \cite{Wa1,Wa2,Wa3}.

In Connes¡¯ program of noncommutative geometry,
for associative unital algebra $\mathcal{A}$ and Hilbert space $\mathcal{H}$ such that there is an algebra
 homomorphism $\pi:\mathcal{A}\rightarrow B(\mathcal{H})$, where $B(\mathcal{H})$
denotes the algebra of bounded operators acting on $\mathcal{H}$,
the role of geometrical
objects is played by spectral triples $(\mathcal{A}; \mathcal{H}; D )$.
An important new feature of such geometries, which is absent in the
commutative case, is the existence of inner fluctuations of the metric.
No doubt it would be extremely interesting to recover scalar curvature, Ricci curvature and other important tensors in both the
 classical setup as well as for the noncommutative geometry or quantum geometry.
Recently, for the metric tensor $g$, Ricci curvature $Ric$ and the scalar curvature $s$,  Dabrowski etc. \cite{DL} recovered the Einstein tensor
\begin{equation}
G:=Ric-\frac{1}{2}s(g)g,
\end{equation}
which directly enters the Einstein field equations with matter, and its dual.
The method of producing these spectral functionals in \cite{DL} is the Wodzicki residue (or be written as  noncommutative residue) density
of a suitable power of the Laplace type operator multiplied by a pair of other differential operators.
Notably, the equivariant Riemannian curvature has the form
 \begin{equation}
R_{g}(X)=(\nabla-\iota(X))^{2}+\mathcal{L}_{X},
\end{equation}
where $\iota(X)$ denotes the contraction operator $\iota(X_{M})$, $\mathcal{L}_{X}=\nabla_{X}+ \mu (X)$ and
the moment map given by  $\mu (X)Y=-\nabla_{Y}X$.
Earlier Jean-Michel Bismut proved the local infinitesimal equivariant index theorem for equivariant
Laplacian  associated with group action \cite{JMB}, and the Bismut Laplacian  given by the formula
 \begin{equation}
H(X)=(D+\frac{1}{4}c(X))^{2}+\mathcal{L}_{X},
\end{equation}
which may be thought of as a quantum analogue of the equivariant Riemannian curvature.
By computing the homologies of the cross-product
algebra $\mathcal{A}(M)\bowtie T$, Dave obtained
an equivariant Connes trace formula as well as extensions of the logarithmic symbols
based on a two cocycle in $H^{2}(S_{log}(M)\times T)$ \cite{SDa}, and constructed an equivariant, delocalized generalization of the
noncommutative residue for manifolds endowed with a group action, which as equivariant noncommutative residue associated to the conjugacy class.
For the equivariant Bismut Laplacian, Wang \cite{Wa4} proved the equivariant Kastler-Kalau-Walze theorems for spin manifolds without boundary
and the equivariant Kastler-Kalau-Walze type theorem with torsion.
And we proved the equivariant Dabrowski-Sitarz-Zalecki type theorems for lower dimensional spin manifolds with (or without) boundary \cite{WW1}.

 Dirac operators with torsion are by now well-established analytical tools in the study of special geometric structures. Ackermann and
  Tolksdorf \cite{AT} proved a generalized version of the well-known Lichnerowicz formula for the square of the
most general Dirac operator with torsion  $D_{T}$ on an even-dimensional spin manifold associated to a metric connection with torsion.
In \cite{PS1}, Pf$\ddot{a}$ffle and Stephan considered orthogonal connections with arbitrary torsion on compact Riemannian manifolds,
and for the induced Dirac operators, twisted Dirac operators and Dirac operators of Chamseddine-Connes type they computed the spectral
action. In \cite{WW2},  we proved the Kastler-Kalau-Walze theorem associated to Dirac operators with torsion on compact manifolds with boundary,
 and gave two kinds of operator-theoretic explanations of the gravitational action in the case of lower dimensional compact manifolds with flat boundary.
Motivated by the spectral Einstein functionals \cite{DL}, equivariant Bismut Laplacian and equivariant Kastler-Kalau-Walze type theorem \cite{Wa4},
 equivariant spectral Einstein functional \cite{WW1},
  we give some new equivariant twisted Kastler-Kalau-Walze type theorems with torsion,
   which are the extension of Kastler-Kalau-Walze type theorems to the equivariant noncommutative with torsion
   realm, and we relate them to the noncommutative residue for lower dimensional manifolds with boundary.

\section{Preliminaries for equivariant twisted Bismut Laplacian with torsion}

The purpose of this section is to specify the equivariant twisted Bismut Laplacian with torsion.
 Let $M$ be a compact oriented Riemannian manifold of even dimension $n$, and let $L$ be a compact Lie group
  acting on $M$. Let  $\mathcal{E}$ be a $L-$equivariant spinor bundle over $M$ with
  $L-$equivariant Hermitian structure and $L-$invariant spin connection $\nabla^{\mathcal{E}}$.
\begin{defn}
A Dirac operator $D$ on a $\mathbb{Z}_{2}$-graded vector bundle $\mathcal{E}$ is a
first-order differential operator of odd parity on $\mathcal{E}$,
 \begin{align}
D:\Gamma(M,\mathcal{E}^{\pm})\rightarrow\Gamma(M,\mathcal{E}^{\mp}),
\end{align}
such that $D^{2}$ is a generalized Laplacian.
\end{defn}
 Let $\nabla^L$ denote the Levi-Civita connection about $g^M$. In the
 fixed orthonormal frame $\{e_1,\cdots,e_n\}$ and natural frames $\{\partial_1,\cdots,\partial_n\}$ of $TM$,
the connection matrix $(\omega_{s,t})$ defined by
\begin{equation}
\nabla^L(e_1,\cdots,e_n)= (e_1,\cdots,e_n)(\omega_{s,t}).
\end{equation}
 Then the Dirac operator has the form
\begin{equation}
D=\sum^n_{i=1}c(e_i)\nabla_{e_i}^{S}=\sum^n_{i=1}c(e_i)\Big[e_i
-\frac{1}{4}\sum_{s,t}\omega_{s,t}(e_i)c(e_s)c(e_t)\Big].
\end{equation}
The Levi-Civita connection
$\nabla^L: \Gamma(TM)\rightarrow \Gamma(T^{*}M\otimes TM)$ on $M$ induces a connection
$\nabla^{S}: \Gamma(S)\rightarrow \Gamma(T^{*}M\otimes S).$ By adding a additional torsion term $t\in\Omega^{1}(M,End TM)$ we
obtain a new covariant derivative
 \begin{equation}
\widetilde{\nabla}:=\nabla^L +t
\end{equation}
on the tangent bundle $TM$. Since $t$ is really a one-form on $M$ with values in the bundle of skew endomorphism $Sk(TM)$ in \cite{GHV},
$\nabla^L$  is in fact compatible with the Riemannian metric $g$ and therefore also induces a connection $\widetilde{\nabla}^{S}:=\nabla^{S}+T$
on the spinor bundle. Here $T\in\Omega^{1}(M, {\rm{ End}} S)$ denotes the `lifted' torsion term $t\in\Omega^{1}(M, {\rm{ End}}(TM))$.
 As Clifford multiplication by any 3-form $T$ is self-adjoint we have
 \begin{equation}
D_{T}\psi=D\psi+\frac{3}{2}T\cdot\psi-\frac{n-1}{2}Y\cdot\psi,
\end{equation}
where the Clifford multiplication by the vector field $Y$ is skew-adjoint,
the hermitian product on the spinor bundle one observes
that $D_{T}$ is symmetric with respect to the natural $L^{2}$-scalar product on spinors if and only
if the vectorial component of the torsion vanishes, $Y\equiv 0$.

Denote by $\mathbb{C}l(M)$ the set of smooth sections of the vector bundle
with fibre over  $x\in M$ the Clifford algebra over the cotangent space of $x$ (a $\mathbb{Z}_{2}$ graded complex algebra
 and $C^{\infty}(M)$-module), we now consider an additional smooth
vector bundle $F$ over $M$ (with $C^{\infty}(M)$-module of smooth sections $W$),  equipped with a connection $\nabla^{F}$,
with corresponding curvature-tensor $R^{F}$. We consider the tensor product vector bundle $S(TM)\otimes F$,
which becomes a Clifford bundle via the definition:
 \begin{equation}
c(a)=\gamma(a)\otimes \texttt{id}_{F},~~~~~a\in \mathbb{C}l(M),
\end{equation}
and which we equip with the compound connection:
 \begin{equation}
\nabla^{ S(TM)\otimes F}= \nabla^{ S(TM)}\otimes \texttt{id}_{ F}+ \texttt{id}_{ S(TM)}\otimes \nabla^{F}.
\end{equation}
The corresponding twisted Dirac operator $D_{F}$ is locally specified as follows:
 \begin{equation}
D_{F}=\sum_{i}^{n}c(e_{i})\nabla^{ S(TM)\otimes F}_{e_{i}}.
\end{equation}
Let $\Phi \in  \Gamma(M, {\rm{ End}} (F))$, motivated by the equivariant Bismut Laplacian and the Connes operator,
we have the following definition.
\begin{defn}
The equivariant twisted Bismut Laplacian with torsion  is the second-order differential operator on $ \Gamma(M,S(TM)\otimes F)$  given by the formula
 \begin{align}
H_{X}=\Big[D_{F}+\frac{1}{4}c(X)+\big(c(T)+c(Y)\big)\otimes \Phi \Big]^{2}+\mathcal{L}^{S(TM)}_{X},
\end{align}
where $\mu^{S(TM)}(X)=\mathcal{L}^{S(TM)}_{X}- \nabla^{S(TM)}_{X}$, $c(X)$ denotes the Clifford action of the dual
of the vector field on $\mathcal{E}$, which satisfies the relation
$
 c(e_{i})c(e_{j})+c(e_{j})c(e_{i})=-2\langle e_{i}, e_{j}  \rangle.
$
When $\Phi=0$, $H_{X}$ is the equivariant twisted Bismut Laplacian;
 and when $F=M\times\mathbb{C} , \Phi={\rm{id}}$, $H_{X}$ is the equivariant Bismut Laplacian with torsion.
\end{defn}
 The following theorem of the equivariant twisted Bismut Laplacian with torsion play a key role in our proof
of  the equivariant twisted Kastler-Kalau-Walze type theorems with torsion.
\begin{thm}
For the the equivariant twisted Bismut Laplacian with torsion $H_{X}$ associated to the orthogonal connection $\nabla^{ S(TM)\otimes F}$
and local natural frames $\{\partial_1,\cdots,\partial_n\}$,
 one has  the following Lichnerowicz formula:
 \begin{align}\label{eq:3.15}
H_{X}
=&-g^{ij}\partial_{i}\partial_{j}
   +g^{ij}\big[X_{j}-2\sigma_{j}^{ S(TM)\otimes F} +\Gamma^{j}
   +c(\partial_{j})\big((c(T)+c(Y)) \otimes \Phi\big)\nonumber\\
&+\Big(\big(c(T)+c(Y)\big)\otimes \Phi\Big) c(\partial_{j})
   +\frac{1}{4}c(\partial_{i})c(X)+\frac{1}{4}c(X)c(\partial_{i})
   \big]\partial_{j}\nonumber\\
    &+g^{ij}\big[-\partial_{i}(\sigma_{j}^{ S(TM)\otimes F} )-\sigma_{i}^{ S(TM)\otimes F}\sigma_{j}^{ S(TM)\otimes F}
      + \Big((c(T)+c(Y))\otimes \Phi \Big)c(\partial_{i})\sigma_{j}^{ S(TM)\otimes F}\nonumber\\
&+ \Gamma^{k}_{ij}\sigma_{k}^{ S(TM)\otimes F}
+c(\partial_{i})\partial_{j}((c(T)+c(Y)) \otimes \Phi)
    +c(\partial_{i})\sigma_{j}^{ S(TM)\otimes F}((c(T)+c(Y)) \otimes \Phi)\nonumber\\
&+\frac{1}{4}c(\partial_{i}) \partial_{j}(c(X))
    +\frac{1}{4}c(\partial_{i})\sigma_{j}^{ S(TM)\otimes F}c(X)
    +\frac{1}{4}c(X)c(\partial_{i})\sigma_{j}^{ S(TM)\otimes F}\big]\nonumber\\
&+\frac{1}{4}(c(T)+c(Y)) \otimes \Phi)c(X)+\frac{1}{4}c(X)(c(T)+c(Y)) \otimes \Phi
   +(c(T)+c(Y))^{2}\otimes \Phi^{2}  \nonumber\\
&+\frac{1}{4}s+\frac{1}{2}\sum_{i\neq j}R^{F}(e_{i},e_{j})c(e_{i})c(e_{j})
    -\frac{1}{16}|X|^{2}+\frac{1}{4}X_{j}\sigma_{j}^{ S(TM)\otimes F}+\mu^{S(TM)}(X),
\end{align}
where $X=\sum_{j=1}^{n}X_{j}\partial_{j} $, $s$ is the scalar curvature and $R^{F}$ denotes the curvature-tensor on $F$.
\end{thm}
\begin{proof}
By Definition 2.2, we have
 \begin{align}
H_{X}=&\Big(D_{F}+(c(T)+c(Y))\otimes \Phi \Big)^{2}
+\frac{1}{4}\Big(D_{F}+(c(T)+c(Y))\otimes \Phi \Big)c(X)\nonumber\\
&+\frac{1}{4}c(X)\Big(D_{F}+(c(T)+c(Y))\otimes \Phi \Big)-\frac{1}{16}|X|^{2}
+\nabla^{S(TM)}_{X}+\mu^{S(TM)}(X).
\end{align}

In what follows, using the Einstein sum convention for repeated index summation:
$\partial^j=g^{ij}\partial_i$,$\sigma^i=g^{ij}\sigma_{j}$, $\Gamma_{k}=g^{ij}\Gamma_{ij}^{k}$,
 and we omit the summation sign, then
\begin{equation}
\nabla_{\partial_i}^{S}=\partial_i+\sigma_i=\partial_i-\frac{1}{4}\sum_{s,t}\omega_{s,t}(e_i)c(e_s)c(e_t) .
\end{equation}
By (2.8) we obtain
 \begin{align}
&\Big(D_{F}+(c(T)+c(Y))\otimes \Phi \Big)^{2}\nonumber\\
=&D_{F}^{2}+D_{F}(c(T)+c(Y))\otimes \Phi+(c(T)+c(Y))\otimes \Phi D_{F}+\big((c(T)+c(Y))\otimes \Phi\big)^{2}\nonumber\\
=&-g^{ij}\partial_{i}\partial_{j}
   +g^{ij}\big[-2\sigma_{j}^{ S(TM)\otimes F} +\Gamma^{j}
   +c(\partial_{i})\big((c(T)+c(Y)) \otimes \Phi\big) +\big((c(T)+c(Y))\otimes \Phi \big)c(\partial_{i})   \big]\partial_{j}\nonumber\\
    &+g^{ij}\big[-\partial_{i}(\sigma_{j}^{ S(TM)\otimes F} )-\sigma_{i}^{ S(TM)\otimes F}\sigma_{j}^{ S(TM)\otimes F}
      + (c(T)+c(Y))\otimes \Phi c(\partial_{i})\sigma_{j}^{ S(TM)\otimes F}\nonumber\\
&+ \Gamma^{k}_{ij}\sigma_{k}^{ S(TM)\otimes F}
+c(\partial_{i})\partial_{j}((c(T)+c(Y)) \otimes \Phi)
    +c(\partial_{i})\sigma_{j}^{ S(TM)\otimes F}((c(T)+c(Y)) \otimes \Phi)\big]\nonumber\\
&+(c(T)+c(Y))^{2}\otimes \Phi^{2}+\frac{1}{4}s+\frac{1}{2}\sum_{i\neq j}R^{F}(e_{i},e_{j})c(e_{i})c(e_{j}).
\end{align}
Also, straightforward computations yield
 \begin{align}
&\frac{1}{4}\Big(D_{F}+(c(T)+c(Y))\otimes \Phi \Big)c(X)+\frac{1}{4}c(X)\Big(D_{F}+(c(T)+c(Y))\otimes \Phi \Big)\nonumber\\
=& g^{ij}\Big(\frac{1}{4}c(\partial_{i})c(X)+\frac{1}{4}c(X)c(\partial_{i})
   \Big)\partial_{j}+\frac{1}{4}(c(T)+c(Y)) \otimes \Phi)c(X)+\frac{1}{4}c(X)(c(T)+c(Y)) \otimes \Phi\nonumber\\
    &+g^{ij}\Big(\frac{1}{4}c(\partial_{i}) \partial_{j}(c(X))
    +\frac{1}{4}c(\partial_{i})\sigma_{j}^{ S(TM)\otimes F}c(X)
    +\frac{1}{4}c(X)c(\partial_{i})\sigma_{j}^{ S(TM)\otimes F}\Big).
\end{align}
Summing up (2.13), (2.14) leads to the desired equality (2.10), and the proof of the Theorem is complete.
\end{proof}

\section{The equivariant twisted Kastler-Kalau-Walze type theorems with torsion}

The aim of this section is to compute the noncommutative residue for  equivariant twisted Bismut Laplacian with torsion
 and get the equivariant twisted Kastler-Kalau-Walze type theorems with torsion on compact manifolds   with  boundary.

\subsection{Noncommutative residue for  equivariant twisted Bismut Laplacian with torsion}

According to the detailed descriptions in \cite{Ac}, we know that the noncommutative residue of a generalized laplacian $\Delta$ is expressed as
 \begin{equation}
(n-2)\Phi_{2}(\Delta)
=(4\pi)^{-\frac{n}{2}}\Gamma(\frac{n}{2})res(\Delta^{-\frac{n}{2}+1}),
 \end{equation}
for $n>2$, where $\Phi_{2}(\Delta)$
denotes the integral over the diagonal part of the second coefficient of the
heat kernel expansion of $\Delta$.
Note that Theorem 1 \cite{KW}, the Wodzicki residue of (certain powers of) an elliptic
differential operator $P$ acting on sections of a complex vector bundle $E$ over a closed compact manifold M was defined by
 \begin{equation}
res(P):=\frac{\Gamma(\frac{n}{2})}{(2\pi)^{\frac{n}{2}}}
\int_{S^{*}M}\mathrm{tr}(\sigma_{-n}^{P}(x,\xi))\mathrm{d}x \mathrm{d}\xi,
\end{equation}
using  $\mathrm{vol}(S^{n-1})^{-1}=\Gamma(\frac{n}{2})$ as a normalizing factor,
then
 \begin{equation}
res(\Delta^{-\frac{n}{2}+1})=\frac{(n-2)}{2}\Phi_{2}(\Delta).
 \end{equation}

Let $V$ be a vector bundle on $M$, any differential operator $P$ of Laplace type has locally the form
 \begin{equation}
P=-\big(g^{ij}\partial_{i}\partial_{j}+A^{i}\partial_{i}+B\big),
\end{equation}
where $\partial_{i}$  is a natural local frame on $TM$ ,
 and $(g^{ij})_{1\leq i,j\leq n}$ is the inverse matrix associated with the metric
matrix  $(g_{ij})_{1\leq i,j\leq n}$ on $M$,
 and $A^{i}$ and $B$ are smooth sections
of $\texttt{End}(V)$ on $M$ \cite{BGV,PBG}. If $P$ is a Laplace type
operator of the form (3.4), then there is a unique
connection $\nabla$ on $V$ and a unique Endomorphism $E$ such that
 \begin{equation}
P=-\Big[g^{ij}(\nabla_{\partial_{i}}\nabla_{\partial_{j}}-
 \nabla_{\nabla^{L}_{\partial_{i}}\partial_{j}})+\widetilde{E}\Big],
\end{equation}
where $\nabla^{L}$ denotes the Levi-Civita connection on $M$. Moreover
(with local frames of $T^{*}M$ and $V$), $\nabla_{\partial_{i}}=\partial_{i}+\omega_{i} $
and $\widetilde{E}$ are related to $g^{ij}$, $A^{i}$, and $B$ through
 \begin{eqnarray}
&&\omega_{i}=\frac{1}{2}g_{ij}\big(A^{i}+g^{kl}\Gamma_{ kl}^{j} \mathrm{Id}\big),\\
&&\widetilde{E}=B-g^{ij}\big(\partial_{i}(\omega_{i})+\omega_{i}\omega_{j}-\omega_{k}\Gamma_{ ij}^{k} \big),
\end{eqnarray}
where $\Gamma_{ kl}^{j}$ is the  Christoffel coefficient of $\nabla^{L}$.
Then for equivariant Connes operators with torsion $H_{X}$, we have the following Lemma.
 \begin{lem}
Let $M$ be $n$-dimensional oriented  compact manifolds  without boundary,
and $H_{X}$ be equivariant twisted Bismut Laplacian with torsion, then
\begin{equation}
\mathrm{Wres}(H_{X}^{-\frac{n-2}{2}})
=\frac{(n-2)\pi^{\frac{n}{2}}}{(\frac{n}{2}-1)!}
\int_{M}\mathrm{tr}(\frac{1}{6}s+\widetilde{E})\mathrm{dvol},
 \end{equation}
where $\mathrm{Wres}$ denote the noncommutative residue, $\mathrm{tr}$ denote trace of operator.
\end{lem}
From Theorem 2.3 and Lemma 3.1, we obtain the main theorem in this section.
\begin{thm}
The noncommutative residue for  equivariant twisted Bismut Laplacian with torsion $H_{X}$ equal to
\begin{align}
\mathrm{Wres}(H_{X}^{-\frac{n-2}{2}})
=&\frac{(n-2)\pi^{\frac{n}{2}}}{(\frac{n}{2}-1)!}\int_{M}\mathrm{dim}(F)2^{n}
 \Big(-\frac{1}{12}s+\frac{1}{2}\mathrm{div}^{g}(X)+g(Y,X){\rm{Tr}}^{F}(\Phi)\nonumber\\
&+\sum _{1\leq\alpha<\beta<\gamma\leq n} 2 T^{2}_{\alpha \beta \gamma} {\rm{Tr}}^{F}(\Phi^{2}) +\frac{1}{2}|Y|^{2}
 {\rm{Tr}}^{F}(\Phi^{2}) \Big)  \mathrm{dvol}_{g},
\end{align}
where $\mathrm{div}^{g}(X)$ denotes the divergence of $X$, and $s$ denotes the scalar curvature.
\end{thm}

 \begin{lem}\label{le:31}
The following identities hold:
 \begin{align}
 &{\rm{Tr}}\big((c(T)+c(Y))c(X)\big)=-g(Y,X){\rm{Tr}}^{ S(TM)}[ {\rm{id}}];\\
& {\rm{Tr}}\big( c(T)c(T)\big)=-\sum _{1\leq\alpha<\beta<\gamma\leq n }T_{\alpha \beta \gamma}^{2 }{\rm{Tr}}^{ S(TM)}[ {\rm{id}}].
\end{align}
\end{lem}

\begin{proof}
Let $T_{\alpha \beta \gamma}=T(e_{\alpha},e_{\beta},e_{\gamma})$,
the first result follows from the  trace properties.
 \begin{align}
 &{\rm{Tr}}\big((c(T)+c(Y))c(X)\big)  \nonumber\\
 =& {\rm{Tr}}(c(T) c(X))(x_{0}) +{\rm{Tr}}(c(Y)c(X))  \nonumber\\
  =& {\rm{Tr}}\Big( \sum _{1\leq\alpha<\beta<\gamma\leq n}
 T_{\alpha \beta \gamma}c(e_{\alpha})c(e_{\beta})c(e_{\gamma})c(X) \Big)
-g(Y,X){\rm{Tr}}^{ S(TM)}[ {\rm{id}}]\nonumber\\
 =& {\rm{Tr}}\Big( \sum _{\substack{1\leq\alpha<\beta<\gamma\leq n \\  l\neq \alpha,\beta,\gamma }}
 T_{\alpha \beta \gamma}c(e_{\alpha})c(e_{\beta})c(e_{\gamma})X_{l}c(e_{l}) \Big) \nonumber\\
&+{\rm{Tr}}\Big( \sum _{\substack{1\leq\alpha<\beta<\gamma\leq n \\  l=\alpha,or~ l=\beta,or~ l=\gamma }}
 T_{\alpha \beta \gamma}c(e_{\alpha})c(e_{\beta})c(e_{\gamma})X_{l}c(e_{l}) \Big)
-g(Y,X){\rm{Tr}}^{ S(TM)}[ {\rm{id}}]\nonumber\\
=&-g(Y,X){\rm{Tr}}^{ S(TM)}[ {\rm{id}}].
\end{align}

In the same way we get
\begin{align}
 {\rm{Tr}}\big( c(T)c(T)\big)
=&{\rm{Tr}}\Big( \sum _{\substack{1\leq\alpha<\beta<\gamma\leq n \\ 1\leq\widetilde{\alpha}<\widetilde{\beta}<\widetilde{\gamma}\leq n}}
 T_{\alpha \beta \gamma}  T_{\widetilde{\alpha} \widetilde{\beta} \widetilde{\gamma}}
c(e_{\alpha })c(e_{\beta}) c(e_{\gamma}) c(e_{\widetilde{\alpha }})c(e_{\widetilde{\beta}}) c(e_{\widetilde{\gamma}})\Big)\nonumber\\
=&{\rm{Tr}}\Big( \sum _{\substack{1\leq\alpha<\beta<\gamma\leq n \\ 1\leq\widetilde{\alpha}<\widetilde{\beta}<\widetilde{\gamma}\leq n}}
 T_{\alpha \beta \gamma}  T_{\widetilde{\alpha} \widetilde{\beta} \widetilde{\gamma}}
c(e_{\beta}) c(e_{\gamma})c(e_{\alpha }) c(e_{\widetilde{\alpha }})c(e_{\widetilde{\beta}}) c(e_{\widetilde{\gamma}})\Big)\nonumber\\
=&{\rm{Tr}}\Big( \sum _{\substack{1\leq\alpha<\beta<\gamma\leq n \\ 1\leq\widetilde{\alpha}<\widetilde{\beta}<\widetilde{\gamma}\leq n}}
 T_{\alpha \beta \gamma}  T_{\widetilde{\alpha} \widetilde{\beta} \widetilde{\gamma}}
c(e_{\beta}) c(e_{\gamma})\big(- c(e_{\widetilde{\alpha }})c(e_{\alpha })-2\delta_{\alpha}^{\widetilde{\alpha}}\big)
c(e_{\widetilde{\beta}}) c(e_{\widetilde{\gamma}})\Big)\nonumber\\
=&-{\rm{Tr}}\Big( \sum _{\substack{1\leq\alpha<\beta<\gamma\leq n \\ 1\leq\widetilde{\alpha}<\widetilde{\beta}<\widetilde{\gamma}\leq n}}
 T_{\alpha \beta \gamma}  T_{\widetilde{\alpha} \widetilde{\beta} \widetilde{\gamma}}
c(e_{\beta}) c(e_{\gamma})c(e_{\widetilde{\alpha }})c(e_{\alpha })
c(e_{\widetilde{\beta}}) c(e_{\widetilde{\gamma}})\Big)\nonumber\\
&-{\rm{Tr}}\Big( \sum _{\substack{1\leq\alpha<\beta<\gamma\leq n \\ 1\leq\widetilde{\alpha}<\widetilde{\beta}<\widetilde{\gamma}\leq n}}
 2 T_{\alpha \beta \gamma}  T_{\widetilde{\alpha} \widetilde{\beta} \widetilde{\gamma}}\delta_{\alpha}^{\widetilde{\alpha}}
c(e_{\beta}) c(e_{\gamma})c(e_{\widetilde{\beta}}) c(e_{\widetilde{\gamma}})\Big)\nonumber\\
=& -{\rm{Tr}}\big( c(T)c(T)\big)-{\rm{Tr}}\Big( \sum _{\substack{1\leq\alpha<\beta<\gamma\leq n \\ 1\leq\widetilde{\alpha}<\widetilde{\beta}<\widetilde{\gamma}\leq n}}
 2 T_{\alpha \beta \gamma}  T_{\widetilde{\alpha} \widetilde{\beta} \widetilde{\gamma}}\delta_{\alpha}^{\widetilde{\alpha}}
c(e_{\beta}) c(e_{\gamma})c(e_{\widetilde{\beta}}) c(e_{\widetilde{\gamma}})\Big)\nonumber\\
&+{\rm{Tr}}\Big( \sum _{\substack{1\leq\alpha<\beta<\gamma\leq n \\ 1\leq\widetilde{\alpha}<\widetilde{\beta}<\widetilde{\gamma}\leq n}}
 2 T_{\alpha \beta \gamma}  T_{\widetilde{\alpha} \widetilde{\beta} \widetilde{\gamma}}\delta_{\alpha}^{\widetilde{\beta }}
c(e_{\beta}) c(e_{\gamma})c(e_{\widetilde{\alpha}}) c(e_{\widetilde{\gamma}})\Big)\nonumber\\
&-{\rm{Tr}}\Big( \sum _{\substack{1\leq\alpha<\beta<\gamma\leq n \\ 1\leq\widetilde{\alpha}<\widetilde{\beta}<\widetilde{\gamma}\leq n}}
 2 T_{\alpha \beta \gamma}  T_{\widetilde{\alpha} \widetilde{\beta} \widetilde{\gamma}}\delta_{\alpha}^{\widetilde{\gamma}}
c(e_{\beta}) c(e_{\gamma})c(e_{\widetilde{\alpha}}) c(e_{\widetilde{\beta}})\Big).
\end{align}
We note that
 \begin{align}\label{eq:4.6}
&{\rm{Tr}}\Big( \sum _{\substack{1\leq\alpha<\beta<\gamma\leq n \\ 1\leq\widetilde{\alpha}<\widetilde{\beta}<\widetilde{\gamma}\leq n}}
T_{\alpha \beta \gamma}  T_{\widetilde{\alpha} \widetilde{\beta} \widetilde{\gamma}}\delta_{\alpha}^{\widetilde{\alpha}}
c(e_{\beta}) c(e_{\gamma})c(e_{\widetilde{\beta}}) c(e_{\widetilde{\gamma}})\Big)\nonumber\\
&=\sum _{\substack{1\leq\alpha<\beta<\gamma\leq n \\ 1\leq\widetilde{\alpha}<\widetilde{\beta}<\widetilde{\gamma}\leq n}}
T_{\alpha \beta \gamma}  T_{\widetilde{\alpha} \widetilde{\beta} \widetilde{\gamma}}\delta_{\alpha}^{\widetilde{\alpha}}
(-\delta_{\beta}^{ \widetilde{\beta}}\delta_{\gamma}^{\widetilde{\gamma}}
+\delta_{\beta}^{\widetilde{\gamma}}\delta_{\gamma}^{\widehat{\beta}}){\rm{Tr}}^{ S(TM)}[ {\rm{id}}]\nonumber\\
&=\sum _{1\leq\alpha<\beta<\gamma\leq n }
T_{\alpha \beta \gamma}^{2 }{\rm{Tr}}^{ S(TM)}[ {\rm{id}}],
\end{align}
 \begin{align}\label{eq:4.6}
&{\rm{Tr}}\Big( \sum _{\substack{1\leq\alpha<\beta<\gamma\leq n \\ 1\leq\widetilde{\alpha}<\widetilde{\beta}<\widetilde{\gamma}\leq n}}
T_{\alpha \beta \gamma}  T_{\widetilde{\alpha} \widetilde{\beta} \widetilde{\gamma}}\delta_{\alpha}^{\widetilde{\beta }}
c(e_{\beta}) c(e_{\gamma})c(e_{\widetilde{\alpha}}) c(e_{\widetilde{\gamma}})\Big)\nonumber\\
&=\sum _{\substack{1\leq\alpha<\beta<\gamma\leq n \\ 1\leq\widetilde{\alpha}<\widetilde{\beta}<\widetilde{\gamma}\leq n}}
T_{\alpha \beta \gamma}  T_{\widetilde{\alpha} \widetilde{\beta} \widetilde{\gamma}}\delta_{\alpha}^{\widetilde{\beta}}
(-\delta_{\beta}^{\widetilde{\alpha} }\delta_{\gamma}^{\widetilde{\gamma}}
+\delta_{\beta}^{\widetilde{\gamma}}\delta_{\gamma}^{\widetilde{\alpha}}){\rm{Tr}}^{ S(TM)}[ {\rm{id}}]\nonumber\\
&=0,
\end{align}
and
 \begin{align}\label{eq:4.6}
&{\rm{Tr}}\Big( \sum _{\substack{1\leq\alpha<\beta<\gamma\leq n \\ 1\leq\widetilde{\alpha}<\widetilde{\beta}<\widetilde{\gamma}\leq n}}
T_{\alpha \beta \gamma}  T_{\widetilde{\alpha} \widetilde{\beta} \widetilde{\gamma}}\delta_{\alpha}^{\widetilde{\gamma}}
c(e_{\beta}) c(e_{\gamma})c(e_{\widetilde{\alpha}}) c(e_{\widetilde{\beta}})\Big)\nonumber\\
&=\sum _{\substack{1\leq\alpha<\beta<\gamma\leq n \\ 1\leq\widetilde{\alpha}<\widetilde{\beta}<\widetilde{\gamma}\leq n}}
T_{\alpha \beta \gamma}  T_{\widetilde{\alpha} \widetilde{\beta} \widetilde{\gamma}}\delta_{\alpha}^{\widetilde{\gamma}}
(-\delta_{\beta}^{\widetilde{\alpha}}\delta_{\gamma}^{\widetilde{\beta}}
+\delta_{\beta}^{\widetilde{\beta}}\delta_{\gamma}^{\widetilde{\alpha}}){\rm{Tr}}^{ S(TM)}[ {\rm{id}}]\nonumber\\
&=0.
\end{align}
Substituting (3.14), (3.15), (3.16) into (3.13), leads to the desired equality (3.11).
\end{proof}

 \begin{lem}\label{le:31}
The following identities hold:
 \begin{align}
&{\rm{Tr}}\Big(\sum_{j}c(e_{j})\nabla_{ e_{j}}^{ S(TM)\otimes F}\big((c(T)+c(Y)) \otimes \Phi\big)\Big)\nonumber\\
=& \Big(\sum _{1\leq\alpha<\beta<\gamma\leq n}(-T_{\alpha j l}+T_{\alpha l j})\langle \nabla _{e_{j}}e_{\alpha}, e_{l}  \rangle
+ \sum _{1\leq l,\beta,j\leq n}T_{l \beta j} \langle \nabla _{e_{j}}e_{\beta}, e_{l}  \rangle\nonumber\\
&-\sum _{1\leq l,j,\gamma\leq n}T_{l j \gamma} \langle \nabla _{e_{j}}e_{\gamma}, e_{l}  \rangle)\Big)
{\rm{Tr}}^{ S(TM)}[ {\rm{Id}}]  {\rm{Tr}}^{F}[ \Phi] -\mathrm{div}^{g}(X){\rm{Tr}}^{ S(TM)}[ {\rm{id}}]{\rm{Tr}}^{F}(\Phi)\nonumber\\
&-{\rm{Tr}}^{ S(TM)}({\rm{Id}})\otimes {\rm{Tr}}^{F}(\nabla_{Y}^{F}(\Phi) )
\end{align}
\end{lem}
\begin{proof}
By the  trace properties and the relation of the Clifford action,  we get
 \begin{align}
 &{\rm{Tr}}\Big(\sum_{j}c(e_{j})\nabla_{ e_{j}}^{ S(TM)\otimes F}\big((c(T)+c(Y)) \otimes \Phi\big)\Big)\nonumber\\
=&{\rm{Tr}}\Big(\sum_{j}c(e_{j}) \big(\nabla_{ e_{j}}^{ S(TM)}(c(T))+c(\nabla_{ e_{j}}^{ L}Y)\big) \otimes \Phi
+ \sum_{j}c(e_{j})(c(T)+c(Y))\otimes\nabla_{ e_{j}}^{F}\Phi\Big)\nonumber\\
=&{\rm{Tr}}\big(\sum_{j}c(e_{j}) \big(\nabla_{ e_{j}}^{ S(TM)}c(T)\big){\rm{Tr}}^{F}(\Phi)
+{\rm{Tr}}\big(\sum_{j}c(e_{j})c(\nabla_{ e_{j}}^{ L}Y)\big) {\rm{Tr}}^{F}(\Phi)
+ {\rm{Tr}}\Big(\sum_{j}c(e_{j})(c(T)+c(Y))\otimes\nabla_{ e_{j}}^{F}\Phi\Big)\nonumber\\
=&{\rm{Tr}}\big(\sum_{j}c(e_{j}) \big(\nabla_{ e_{j}}^{ S(TM)}c(T)\big){\rm{Tr}}^{F}(\Phi)
-\mathrm{div}^{g}(X){\rm{Tr}}^{ S(TM)}[ {\rm{id}}]{\rm{Tr}}^{F}(\Phi)
+ {\rm{Tr}}\Big(\sum_{j}c(e_{j})(c(T)+c(Y))\otimes\nabla_{ e_{j}}^{F}\Phi\Big)\nonumber\\
=&{\rm{Tr}}\big(\sum_{j}c(e_{j}) \big(\nabla_{ e_{j}}^{ S(TM)}c(T)\big){\rm{Tr}}^{F}(\Phi)
-\mathrm{div}^{g}(X){\rm{Tr}}^{ S(TM)}[ {\rm{id}}]{\rm{Tr}}^{F}(\Phi)
-{\rm{Tr}}({\rm{Id}})\otimes {\rm{Tr}}(\nabla_{Y}^{F}(\Phi) ).
\end{align}
And a simple computation shows
 \begin{align}
&{\rm{Tr}}\Big(\sum_{j}c(e_{j}) \big(\nabla_{ e_{j}}^{ S(TM)}c(T)\big)\Big) \nonumber\\
=&{\rm{Tr}}\Big( \sum _{1\leq\alpha<\beta<\gamma\leq n}c(e_{j})
\big( e_{j}(T_{\alpha \beta \gamma})c(e_{\alpha})c(e_{\beta})c(e_{\gamma})
 +T_{\alpha \beta \gamma} c(\nabla_{ e_{j}}^{L}e_{\alpha} )c(e_{\alpha}) c(e_{\gamma}) \nonumber\\
 &+T_{\alpha \beta \gamma}c(e_{\alpha})c(\nabla_{ e_{j}}^{L}e_{\beta} ) c(e_{\gamma})
 +T_{\alpha \beta \gamma}c(e_{\alpha})c(e_{\beta})c(\nabla_{ e_{j}}^{L}e_{\gamma} )
\big)\Big)\nonumber\\
=&{\rm{Tr}}\Big( \sum _{1\leq\alpha<\beta<\gamma\leq n}
\big( T_{\alpha \beta \gamma}c(e_{j}) c(\nabla_{ e_{j}}^{L}e_{\alpha} )c(e_{\beta}) c(e_{\gamma}) \nonumber\\
 &+T_{\alpha \beta \gamma}c(e_{j})c(e_{\alpha})c(\nabla_{ e_{j}}^{L}e_{\beta} ) c(e_{\gamma})
 +T_{\alpha \beta \gamma}c(e_{j})c(e_{\alpha})c(e_{\beta})c(\nabla_{ e_{j}}^{L}e_{\gamma} )
\big)\Big)\nonumber\\
=&{\rm{Tr}}\Big( \sum _{1\leq\alpha<\beta<\gamma\leq n}
\big( T_{\alpha \beta \gamma} \sum _{l}\langle \nabla _{e_{j}}e_{\alpha}, e_{l}  \rangle
c(e_{j}) c(e_{l} )c(e_{\beta}) c(e_{\gamma}) \nonumber\\
 &+T_{\alpha \beta \gamma} \sum _{l}\langle \nabla _{e_{j}}e_{\beta}, e_{l}  \rangle
  c(e_{j})c(e_{\alpha})c( e_{l} ) c(e_{\gamma})
 +T_{\alpha \beta \gamma}\sum _{l}\langle \nabla _{e_{j}}e_{\gamma}, e_{l}  \rangle
 c(e_{j})c(e_{\alpha})c(e_{\beta})c(e_{l} )
\big)\Big).
\end{align}
Also, straightforward computations yield
 \begin{align}
&{\rm{Tr}}\Big( \sum _{1\leq\alpha<\beta<\gamma\leq n}
 T_{\alpha \beta \gamma} \sum _{l}\langle \nabla _{e_{j}}e_{\alpha}, e_{l}  \rangle
c(e_{j}) c(e_{l} )c(e_{\beta}) c(e_{\gamma}) \Big)\nonumber\\
=&\sum _{1\leq\alpha<\beta<\gamma\leq n}
 T_{\alpha \beta \gamma} \sum _{l}\langle \nabla _{e_{j}}e_{\alpha}, e_{l}  \rangle
(-\delta_{j}^{ \beta}\delta_{l}^{\gamma}+\delta_{j}^{\gamma}\delta_{l}^{\beta}){\rm{Tr}}^{ S(TM)}[ {\rm{id}}]\nonumber\\
=&\Big(\sum _{1\leq\alpha<\beta<\gamma\leq n}
 -T_{\alpha \beta \gamma} \sum _{l}\langle \nabla _{e_{j}}e_{\alpha}, e_{l}  \rangle
\delta_{j}^{ \beta}\delta_{l}^{\gamma}
+\sum _{1\leq\alpha<\beta<\gamma\leq n}
 T_{\alpha \beta \gamma} \sum _{l}\langle \nabla _{e_{j}}e_{\alpha}, e_{l}  \rangle
\delta_{j}^{\gamma}\delta_{l}^{\beta}
\Big){\rm{Tr}}^{ S(TM)}[ {\rm{id}}]
\nonumber\\
=& \sum _{1\leq\alpha<j<l\leq n}(-T_{\alpha j l}+T_{\alpha l j})\langle \nabla _{e_{j}}e_{\alpha}, e_{l}  \rangle
{\rm{Tr}}^{ S(TM)}[ {\rm{id}}]=\sum _{1\leq\alpha<j<l\leq n}(-2T_{\alpha j l})\langle \nabla _{e_{j}}e_{\alpha}, e_{l}  \rangle
{\rm{Tr}}^{ S(TM)}[ {\rm{id}}],
\end{align}
and
 \begin{align}
&{\rm{Tr}}\Big( \sum _{1\leq\alpha<\beta<\gamma\leq n}
 T_{\alpha \beta \gamma} \sum _{l}\langle \nabla _{e_{j}}e_{\beta}, e_{l}  \rangle
c(e_{j}) c(e_{\alpha } )c(e_{l}) c(e_{\gamma}) \Big)\nonumber\\
=&\sum _{1\leq\alpha<\beta<\gamma\leq n}
 T_{\alpha \beta \gamma} \sum _{l}\langle \nabla _{e_{j}}e_{\beta}, e_{l}  \rangle
(-\delta_{j}^{l}\delta_{\alpha}^{\gamma}+\delta_{j}^{\gamma}\delta_{\alpha}^{l}){\rm{Tr}}^{ S(TM)}[ {\rm{id}}]\nonumber\\
=&\Big(\sum _{1\leq\alpha<\beta<\gamma\leq n}
 -T_{\alpha \beta \gamma} \sum _{l}\langle \nabla _{e_{j}}e_{\beta}, e_{l}  \rangle
\delta_{j}^{l}\delta_{\alpha}^{\gamma}
+\sum _{1\leq\alpha<\beta<\gamma\leq n}
 T_{\alpha \beta \gamma} \sum _{l}\langle \nabla _{e_{j}}e_{\alpha}, e_{l}  \rangle
\delta_{j}^{\gamma}\delta_{\alpha}^{l}
\Big){\rm{Tr}}^{ S(TM)}[ {\rm{id}}]
\nonumber\\
=& \sum _{1\leq l,\beta,j\leq n}T_{l \beta j} \langle \nabla _{e_{j}}e_{\beta}, e_{l}  \rangle
{\rm{Tr}}^{ S(TM)}[ {\rm{id}}].
\end{align}
Similarly,
 \begin{align}
&{\rm{Tr}}\Big( \sum _{1\leq\alpha<\beta<\gamma\leq n}
 T_{\alpha \beta \gamma} \sum _{l}\langle \nabla _{e_{j}}e_{\gamma}, e_{l}  \rangle
c(e_{j}) c(e_{\alpha} )c(e_{\beta}) c(e_{l}) \Big)
 =-\sum _{1\leq l,j,\gamma\leq n}T_{l j \gamma} \langle \nabla _{e_{j}}e_{\gamma}, e_{l}  \rangle
{\rm{Tr}}^{ S(TM)}[ {\rm{id}}].
\end{align}
Summing up (3.18)- (3.22) leads to (3.17), and the proof of the Lemma is complete.
\end{proof}
We shall prove Theorem 3.2 at the end of this section.

 {\bf Proof of Theorem 3.2}
\begin{proof}
Let
 \begin{align}
H_{X}=\Big[D_{F}+\frac{1}{4}c(X)+\big(c(T)+c(Y)\big)\otimes \Phi \Big]^{2}+\mathcal{L}^{S(TM)}_{X}=
-\big(g^{ij}\partial_{i}\partial_{j}+A^{i}\partial_{i}+B\big).
\end{align}
From Theorem 2.3 we have
 \begin{align}
A^{i}=&-g^{ij}\big[X_{j}-2\sigma_{j}^{ S(TM)\otimes F} +\Gamma^{j}
   +c(\partial_{j})\big((c(T)+c(Y)) \otimes \Phi\big)\nonumber\\
&+\big(c(T)+c(Y)\big)\otimes \Phi c(\partial_{j})
   +\frac{1}{4}c(\partial_{j})c(X)+\frac{1}{4}c(X)c(\partial_{j})
   \big],
\end{align}
 \begin{align}
B=&g^{ij}\big[-\partial_{i}(\sigma_{j}^{ S(TM)\otimes F} )-\sigma_{i}^{ S(TM)\otimes F}\sigma_{j}^{ S(TM)\otimes F}
      + (c(T)+c(Y))\otimes \Phi c(\partial_{i})\sigma_{j}^{ S(TM)\otimes F}\nonumber\\
&+ \Gamma^{k}_{ij}\sigma_{k}^{ S(TM)\otimes F}
+c(\partial_{i})\partial_{j}((c(T)+c(Y)) \otimes \Phi)
    +c(\partial_{j})\sigma_{j}^{ S(TM)\otimes F}((c(T)+c(Y)) \otimes \Phi)\nonumber\\
&+\frac{1}{4}c(\partial_{i}) \partial_{j}(c(X))
    +\frac{1}{4}c(\partial_{i})\sigma_{j}^{ S(TM)\otimes F}c(X)
    +\frac{1}{4}c(X)c(\partial_{i})\sigma_{j}^{ S(TM)\otimes F}\big]\nonumber\\
&-\frac{1}{4}(c(T)+c(Y)) \otimes \Phi)c(X)-\frac{1}{4}c(X)(c(T)+c(Y)) \otimes \Phi
   -(c(T)+c(Y))^{2}\otimes \Phi^{2}  \nonumber\\
&-\frac{1}{4}s-\frac{1}{2}\sum_{i\neq j}R^{F}(e_{i},e_{j})c(e_{i})c(e_{j})
    +\frac{1}{16}|X|^{2}-\frac{1}{4}X_{j}\sigma_{j}^{ S(TM)\otimes F}-\mu^{S(TM)}(X),
\end{align}
and
 \begin{align}
\omega_{i}=\frac{1}{2}g_{ij}\big(A^{i}+g^{kl}\Gamma_{ kl}^{j} \mathrm{Id}\big).
\end{align}
It follows from (3.7) that $\widetilde{E}=B-g^{ij}\big(\partial_{i}(\omega_{i})+\omega_{i}\omega_{j}-\omega_{k}\Gamma_{ ij}^{k} \big)$.
In view of (3.6), it suffices to show that
 \begin{align}\label{eq:3.12}
g^{ij}\partial_{i}(\omega_{i})&= \partial_{i}\Big[ \frac{1}{2}\sum_{i,j}g_{ij}\Big( 2\sigma_{j}^{ S(TM)\otimes F}-\Gamma^{j}
-c(\partial_{i})\big(c(T)+c(Y)\big)\otimes \Phi
- \big(c(T)+c(Y)\big)\otimes \Phi c(\partial_{i}) +g^{kl}\Gamma^{j}_{kl}\mathrm{id} \Big) \Big].
\end{align}
Also, straightforward computations yield
 \begin{align}\label{eq:3.12}
\omega_{i}\omega_{j}=&\frac{1}{2}g_{ik}\big(A^{i}+g^{kl}\Gamma_{ kl}^{j} \mathrm{Id}\big)\frac{1}{2}g_{jl}\big(A^{j}
+g^{kl}\Gamma_{ kl}^{j}\big)\nonumber\\
=&\frac{1}{2}g_{ik}\Big[ g^{ij}\big[X_{i}-2\sigma_{i}^{ S(TM)\otimes F} -\Gamma^{i}
   +c(\partial_{i})\big((c(T)+c(Y)) \otimes \Phi\big)\nonumber\\
&+\big(c(T)+c(Y)\big)\otimes \Phi c(\partial_{i})
   +\frac{1}{4}c(\partial_{i})c(X)+\frac{1}{4}c(X)c(\partial_{i})
   \big] +g^{kl}\Gamma_{ kl}^{i} \mathrm{Id} \Big]\nonumber\\
&\times \frac{1}{2}g_{jl} \Big[g^{ij}\big[X_{j}-2\sigma_{j}^{ S(TM)\otimes F}-\Gamma^{j}
   +c(\partial_{j})\big((c(T)+c(Y)) \otimes \Phi\big)\nonumber\\
&+\big(c(T)+c(Y)\big)\otimes \Phi c(\partial_{j})
   +\frac{1}{4}c(\partial_{j})c(X)+\frac{1}{4}c(X)c(\partial_{j})
   \big]+g^{kl}\Gamma_{ kl}^{j} \Big]\nonumber\\
=&\frac{1}{4}g_{ik}g_{jl}X^{i}X^{j}
   -\frac{1}{4}g_{ik}g_{jl}X^{j}\Big(  2\sigma^{j}_{ S(TM)\otimes F}-c(\partial_{j})\big((c(T)+c(Y)) \otimes \Phi\big)\nonumber\\
     &-\big(c(T)+c(Y)\big)\otimes \Phi c(\partial_{j})
   -\frac{1}{4}c(\partial_{j})c(X)-\frac{1}{4}c(X)c(\partial_{j}) \Big)\nonumber\\
&+\frac{1}{4}g_{ik}g_{jl}\sigma^{j}_{ S(TM)\otimes F}\Big(-X^{i}+2\sigma^{i}_{ S(TM)\otimes F}
-c(\partial_{j})\big((c(T)+c(Y)) \otimes \Phi\big)\nonumber\\
     &-\big(c(T)+c(Y)\big)\otimes \Phi c(\partial_{j})
   -\frac{1}{4}c(\partial_{j})c(X)-\frac{1}{4}c(X)c(\partial_{j}) \Big)\nonumber\\
&+\frac{1}{4}g_{ik}g_{jl} \Big(c(\partial_{j})\big((c(T)+c(Y)) \otimes \Phi\big)+\big(c(T)+c(Y)\big)\otimes \Phi c(\partial_{j})
      \Big)\nonumber\\
&\times\Big(-X^{i}+2\sigma^{i}_{ S(TM)\otimes F}
-c(\partial_{j})\big((c(T)+c(Y)) \otimes \Phi\big)\nonumber\\
     &-\big(c(T)+c(Y)\big)\otimes \Phi c(\partial_{j})
   -\frac{1}{4}c(\partial_{j})c(X)-\frac{1}{4}c(X)c(\partial_{j}) \Big)\nonumber\\
    &-\frac{1}{16}g_{ik}g_{jl}\Big( c(\partial_{j})c(X)+c(X)c(\partial_{j})  \Big)
 \Big(-X^{i}+2\sigma^{i}_{ S(TM)\otimes F}
   -c(\partial_{j})\big((c(T)+c(Y)) \otimes \Phi\big)\nonumber\\
     &-\big(c(T)+c(Y)\big)\otimes \Phi c(\partial_{j})
   -\frac{1}{4}c(\partial_{j})c(X)-\frac{1}{4}c(X)c(\partial_{j}) \Big).
         \end{align}
Combining (3.25), (3.27) and (3.28), we obtain

 \begin{align}\label{eq:3.18}
E=&-\frac{1}{4}s-\frac{1}{2}\sum_{i\neq j}R^{F}(e_{i},e_{j})c(e_{i})c(e_{j})
    +\frac{1}{16}|X|^{2}-\frac{1}{4}X_{j}\sigma_{j}^{ S(TM)\otimes F}-\mu^{S(TM)}(X)\nonumber\\
&+g^{ij}\big[-\partial_{i}(\sigma_{j}^{ S(TM)\otimes F} )-\sigma_{i}^{ S(TM)\otimes F}\sigma_{j}^{ S(TM)\otimes F}
      + (c(T)+c(Y))\otimes \Phi c(\partial_{i})\sigma_{j}^{ S(TM)\otimes F}\nonumber\\
&+ \Gamma^{k}_{ij}\sigma_{k}^{ S(TM)\otimes F}
+c(\partial_{i})\partial_{j}((c(T)+c(Y)) \otimes \Phi)
    +c(\partial_{j})\sigma_{j}^{ S(TM)\otimes F}((c(T)+c(Y)) \otimes \Phi)\nonumber\\
&+\frac{1}{4}c(\partial_{i}) \partial_{j}(c(X))
    +\frac{1}{4}c(\partial_{i})\sigma_{j}^{ S(TM)\otimes F}c(X)
    +\frac{1}{4}c(X)c(\partial_{i})\sigma_{j}^{ S(TM)\otimes F}\big]\nonumber\\
&-\frac{1}{4}(c(T)+c(Y)) \otimes \Phi)c(X)-\frac{1}{4}c(X)(c(T)+c(Y)) \otimes \Phi
   -(c(T)+c(Y))^{2}\otimes \Phi^{2}  \nonumber\\
&-\partial_{i}\Big[ \frac{1}{2}\sum_{i,j}g_{ij}\Big( 2\sigma_{j}^{ S(TM)\otimes F}-\Gamma^{j}
-c(\partial_{i})\big(c(T)+c(Y)\big)\otimes \Phi
- \big(c(T)+c(Y)\big)\otimes \Phi c(\partial_{i}) +g^{kl}\Gamma^{j}_{kl}\mathrm{id} \Big) \Big]\nonumber\\
&-\frac{1}{4}g_{ij}g_{ik}g_{jl}X^{i}X^{j}
   -\frac{1}{4}g_{ij}g_{ik}g_{jl}X^{j}\Big(  2\sigma^{j}_{ S(TM)\otimes F}-c(\partial_{j})\big((c(T)+c(Y)) \otimes \Phi\big)\nonumber\\
&-\big(c(T)+c(Y)\big)\otimes \Phi c(\partial_{j})
   -\frac{1}{4}c(\partial_{j})c(X)-\frac{1}{4}c(X)c(\partial_{j}) \Big)\nonumber\\
&-\frac{1}{4}g_{ij}g_{ik}g_{jl}\sigma^{j}_{ S(TM)\otimes F}\Big(-X^{i}+2\sigma^{i}_{ S(TM)\otimes F}
-c(\partial_{j})\big((c(T)+c(Y)) \otimes \Phi\big)\nonumber\\
&-\big(c(T)+c(Y)\big)\otimes \Phi c(\partial_{j})
   -\frac{1}{4}c(\partial_{j})c(X)-\frac{1}{4}c(X)c(\partial_{j}) \Big)\nonumber\\
&-\frac{1}{4}g_{ij}g_{ik}g_{jl} \Big(c(\partial_{j})\big((c(T)+c(Y)) \otimes \Phi\big)+\big(c(T)+c(Y)\big)\otimes \Phi c(\partial_{j})
      \Big)\nonumber\\
&\times\Big(-X^{i}+2\sigma^{i}_{ S(TM)\otimes F}
-c(\partial_{j})\big((c(T)+c(Y)) \otimes \Phi\big)\nonumber\\
&-\big(c(T)+c(Y)\big)\otimes \Phi c(\partial_{j})
   -\frac{1}{4}c(\partial_{j})c(X)-\frac{1}{4}c(X)c(\partial_{j}) \Big)\nonumber\\
&+\frac{1}{16}g_{ij}g_{ik}g_{jl}\Big( c(\partial_{j})c(X)+c(X)c(\partial_{j})  \Big)
 \Big(-X^{i}+2\sigma^{i}_{ S(TM)\otimes F} -c(\partial_{j})\big((c(T)+c(Y)) \otimes \Phi\big)\nonumber\\
 &-\big(c(T)+c(Y)\big)\otimes \Phi c(\partial_{j})
   -\frac{1}{4}c(\partial_{j})c(X)-\frac{1}{4}c(X)c(\partial_{j}) \Big)\nonumber\\
 &+\frac{1}{2}\sum_{i,j}g_{ij}\Big[\Big( 2\sigma_{j}^{ S(TM)\otimes F} )-\Gamma^{j}
-c(\partial_{i})\big(c(T)+c(Y)\big)\otimes \Phi
- \big(c(T)+c(Y)\big)\otimes \Phi c(\partial_{i}) +g^{kl}\Gamma^{j}_{kl}\mathrm{id} \Big]  \Gamma^{k}_{ij}.
   \end{align}

Since $E$ is globally defined on $M$, so we can perform
computations of $E$ in normal coordinates. In terms of normal coordinates about $x_{0}$ one has:
$\sigma^{j}_{S(TM)}(x_{0})=0$ $e_{j}\big(c(e_{i})\big)(x_{0})=0$, $\Gamma^{k}(x_{0})=0$, we conclude that
\begin{align}\label{eq:3.19}
E(x_{0})=&-\frac{1}{4}s-\frac{1}{2}\sum_{i\neq j}R^{F}(e_{i},e_{j})c(e_{i})c(e_{j})
         +\frac{1}{4}\mathrm{div}^{g}(X)-\mu^{S(TM)}(X)(x_{0})\nonumber\\
    &-\frac{1}{4}\sum_{j}c(e_{j})c(\nabla_{ e_{j}}^{L}X )(x_{0})
    -\sum_{j}c(e_{j})\nabla_{ e_{j}}^{ S(TM)\otimes F}\big((c(T)+c(Y)) \otimes \Phi\big)\nonumber\\
    &-\frac{1}{2}\big((c(T)+c(Y))\otimes \Phi \big)c(X)-\frac{1}{2}c(X)\big(c(T)+c(Y)\big)\otimes \Phi \nonumber\\
    &+\frac{1}{2}\sum_{j}\nabla_{ e_{j}}^{ S(TM)\otimes F}\big((c(T)+c(Y)) \otimes \Phi \big)c(e_{j})
    +\frac{1}{2}\sum_{j}c(e_{j}) \otimes \nabla_{ e_{j}}^{ S(TM)\otimes F} \big((c(T)+c(Y))\otimes \Phi)\big) \nonumber\\
    &-\frac{1}{4}\sum_{j}\big((c(T)+c(Y)) c(e_{j})
    +c(e_{j})(c(T)+c(Y))\big)^2\otimes \Phi^{2}- (c(T)+c(Y))^2\otimes \Phi^2.
\end{align}
Next we compute ${\rm{Tr}}(E(x_{0}))$.
A simple computation shows
\begin{align}
\mathrm{Tr}^{S(TM)} \Big[e_{j}\Big(\frac{1}{4}c(X)\Big)c(e_{j})-c(e_{j})e_{j}\Big(\frac{1}{4}c(X)\Big)
\Big](x_{0})=0,
\end{align}
and
 \begin{align}
-g^{ij} \partial_{i}\big(  \frac{1}{2}g_{jk}X_{k}\big)(x_{0})
=&-\frac{1}{2}\sum_{j=1}^{n}\partial_{j}\big( g_{jk}X_{k}\big)(x_{0})\nonumber\\
=&-\frac{1}{2}\sum_{j=1}^{n}\Big(\partial_{j}\big( g_{jk}\big)(x_{0})X_{k}+ g_{jk}(x_{0})\partial_{j}\big(X_{k}\big)(x_{0})   \Big)\nonumber\\
=&-\frac{1}{2}\sum_{j=1}^{n}\partial_{j}(X_{j})(x_{0})=-\frac{1}{2}\mathrm{div}(X).
\end{align}
Similarly,
 \begin{align}
\sum_{j=1}^{n}\langle \nabla_{ e_{j}}^{TM}X, e_{j} \rangle(x_{0})
=&\sum_{j,k=1}^{n}\langle \nabla_{ e_{j}}^{TM}(X_{k}\partial_{k}), e_{j} \rangle(x_{0})\nonumber\\
=&\sum_{j,k=1}^{n}\langle e_{j}(X_{k})\partial_{k}+X_{k}\nabla_{ e_{j}}^{TM}(\partial_{k}), e_{j} \rangle(x_{0})\nonumber\\
=&\sum_{j=1}^{n}\partial_{j}(X_{j})(x_{0})= \mathrm{div}(X),
\end{align}
and
 \begin{align}
-\frac{1}{4}g^{ij}g_{ik}g_{jl}X_{k}X_{l}(x_{0}) =&-\sum_{i,j,k,l}\frac{1}{4}g^{ij}g_{ik}g_{jl}X_{k}X_{l}(x_{0})\nonumber\\
=&-\sum_{ j,k,l}\frac{1}{4}\delta^{jk}g_{jl}X_{k}X_{l}(x_{0})\nonumber\\
=&-\frac{1}{4}\sum_{ j }X_{j}^{2}(x_{0})=-\frac{1}{4}|X|^{2}.
\end{align}
By direct computations, we obtain
 \begin{align}
{\rm{Tr}}(c(T)c( \partial_{n}))(x_{0})=
&{\rm{Tr}}\Big( \sum _{1\leq j<k<l\leq n}  T(e_{j},e_{k},e_{l})c(e_{j})c(e_{k})c(e_{l})c( \partial_{n})\Big)\nonumber\\
=&{\rm{Tr}}\Big( \sum _{1\leq j<k<l< n}  T(e_{j},e_{k},e_{l})c(e_{j})c(e_{k})c(e_{l})c( \partial_{n})\Big)\nonumber\\
 &+{\rm{Tr}}\Big( \sum _{1\leq j<k<l=n}  T(e_{j},e_{k},e_{l})c(e_{j})c(e_{k})c(e_{n})c( \partial_{n})\Big)=0,
\end{align}
and
 \begin{align}
{\rm{Tr}}(c(Y)c( \partial_{n}))(x_{0})=-g(Y, \partial_{n}){\rm{Tr}}^{ S(TM)\otimes F}[ {\rm{id}}],
\end{align}
 \begin{align}
 {\rm{Tr}}(\sum_{j}c(e_{j})c(\nabla_{ e_{j}}^{L}X ))(x_{0})=-{\rm{Tr}}(\sum_{j}g(e_{j},\nabla_{ e_{j}}^{L}X ))(x_{0})
=-\mathrm{div}^{g}(X){\rm{Tr}}^{ S(TM)}[ {\rm{id}}].
\end{align}

 \begin{align}
{\rm{Tr}}\big(\mu^{S(TM)}(X)\big)(x_{0})=& {\rm{Tr}}(\mathcal{L}^{S(TM)}_{X}- \nabla^{S(TM)}_{X})(x_{0})
={\rm{Tr}}(-\frac{1}{4} c(dX^{*}))(x_{0})\nonumber\\
=&{\rm{Tr}}(-\frac{1}{4} \sum _{j<l}dX^{*}(e_{j},e_{l} )c(e_{j})c(e_{l}) )(x_{0})=0,
\end{align}
where $X^{*}=g(X,\cdot)$.
By Lemma 3.3 and Lemma 3.4, summing up (3.31)-(3.38) we obtain
 \begin{align}\label{eq:4.6}
{\rm{Tr}}(E(x_{0}))=&2^{n}\mathrm{dim}F
 \Big(-\frac{1}{4}s+\frac{1}{2}\mathrm{div}^{g}(X)+g(Y,X){\rm{Tr}}^{ F}(\Phi)\nonumber\\
&+\sum _{1\leq\alpha<\beta<\gamma\leq n} 2 T^{2}_{\alpha \beta \gamma} {\rm{Tr}}^{ F}(\Phi^{2})
+\frac{1}{2}|Y|^{2} {\rm{Tr}}^{ F}(\Phi^{2}) \Big).
\end{align}
Combining (3.8) and (3.39) leads to the desired equality, and the proof of the Theorem 3.2 is complete.
\end{proof}

\subsection{The equivariant twisted Kastler-Kalau-Walze type theorems with torsion on compact manifolds   with  boundary }
In this section, for compact oriented manifold $M$ with boundary $\partial M$,
to define the noncommutative residue associated with equivariant twisted Bismut Laplacian  with torsion $H_{X}$,
some basic facts and formulae about Boutet de Monvel's calculus are recalled as follows.

Let $$ F:L^2({\bf R}_t)\rightarrow L^2({\bf R}_v);~F(u)(v)=\int e^{-ivt}u(t)\texttt{d}t$$ denote the Fourier transformation and
$\varphi(\overline{{\bf R}^+}) =r^+\varphi({\bf R})$ (similarly define $\varphi(\overline{{\bf R}^-}$)), where $\varphi({\bf R})$
denotes the Schwartz space and
  \begin{equation}
r^{+}:C^\infty ({\bf R})\rightarrow C^\infty (\overline{{\bf R}^+});~ f\rightarrow f|\overline{{\bf R}^+};~
 \overline{{\bf R}^+}=\{x\geq0;x\in {\bf R}\}.
\end{equation}
We define $H^+=F(\varphi(\overline{{\bf R}^+}));~ H^-_0=F(\varphi(\overline{{\bf R}^-}))$ which are orthogonal to each other. We have the following
 property: $h\in H^+~(H^-_0)$ iff $h\in C^\infty({\bf R})$ which has an analytic extension to the lower (upper) complex
half-plane $\{{\rm Im}\xi<0\}~(\{{\rm Im}\xi>0\})$ such that for all nonnegative integer $l$,
 \begin{equation}
\frac{\texttt{d}^{l}h}{\texttt{d}\xi^l}(\xi)\sim\sum^{\infty}_{k=1}\frac{\texttt{d}^l}{\texttt{d}\xi^l}(\frac{c_k}{\xi^k})
\end{equation}
as $|\xi|\rightarrow +\infty,{\rm Im}\xi\leq0~({\rm Im}\xi\geq0)$.

 Let $H'$ be the space of all polynomials and $H^-=H^-_0\bigoplus H';~H=H^+\bigoplus H^-.$ Denote by $\pi^+~(\pi^-)$ respectively the
 projection on $H^+~(H^-)$. For calculations, we take $H=\widetilde H=\{$rational functions having no poles on the real axis$\}$ ($\tilde{H}$
 is a dense set in the topology of $H$). Then on $\tilde{H}$,
 \begin{equation}
\pi^+h(\xi_0)=\frac{1}{2\pi i}\lim_{u\rightarrow 0^{-}}\int_{\Gamma^+}\frac{h(\xi)}{\xi_0+iu-\xi}\texttt{d}\xi,
\end{equation}
where $\Gamma^+$ is a Jordan close curve included ${\rm Im}\xi>0$ surrounding all the singularities of $h$ in the upper half-plane and
$\xi_0\in {\bf R}$. Similarly, define $\pi^{'}$ on $\tilde{H}$,
 \begin{equation}
\pi'h=\frac{1}{2\pi}\int_{\Gamma^+}h(\xi)\texttt{d}\xi.
\end{equation}
So, $\pi'(H^-)=0$. For $h\in H\bigcap L^1(R)$, $\pi'h=\frac{1}{2\pi}\int_{R}h(v)\texttt{d}v$ and for $h\in H^+\bigcap L^1(R)$, $\pi'h=0$.
Denote by $\mathcal{B}$ Boutet de Monvel's algebra (for details, see Section 2 of \cite{Wa1}).

An operator of order $m\in {\bf Z}$ and type $d$ is a matrix
$$A=\left(\begin{array}{lcr}
  \pi^+P+G  & K  \\
   T  &  S
\end{array}\right):
\begin{array}{cc}
\   C^{\infty}(X,E_1)\\
 \   \bigoplus\\
 \   C^{\infty}(\partial{X},F_1)
\end{array}
\longrightarrow
\begin{array}{cc}
\   C^{\infty}(X,E_2)\\
\   \bigoplus\\
 \   C^{\infty}(\partial{X},F_2)
\end{array}.
$$
where $X$ is a manifold with boundary $\partial X$ and
$E_1,E_2~(F_1,F_2)$ are vector bundles over $X~(\partial X
)$.~Here,~$P:C^{\infty}_0(\Omega,\overline {E_1})\rightarrow
C^{\infty}(\Omega,\overline {E_2})$ is a classical
pseudodifferential operator of order $m$ on $\Omega$, where
$\Omega$ is an open neighborhood of $X$ and
$\overline{E_i}|X=E_i~(i=1,2)$. $P$ has an extension:
$~{\cal{E'}}(\Omega,\overline {E_1})\rightarrow
{\cal{D'}}(\Omega,\overline {E_2})$, where
${\cal{E'}}(\Omega,\overline {E_1})~({\cal{D'}}(\Omega,\overline
{E_2}))$ is the dual space of $C^{\infty}(\Omega,\overline
{E_1})~(C^{\infty}_0(\Omega,\overline {E_2}))$. Let
$e^+:C^{\infty}(X,{E_1})\rightarrow{\cal{E'}}(\Omega,\overline
{E_1})$ denote extension by zero from $X$ to $\Omega$ and
$r^+:{\cal{D'}}(\Omega,\overline{E_2})\rightarrow
{\cal{D'}}(\Omega, {E_2})$ denote the restriction from $\Omega$ to
$X$, then define
$$\pi^+P=r^+Pe^+:C^{\infty}(X,{E_1})\rightarrow {\cal{D'}}(\Omega,
{E_2}).$$
In addition, $P$ is supposed to have the
transmission property; this means that, for all $j,k,\alpha$, the
homogeneous component $p_j$ of order $j$ in the asymptotic
expansion of the
symbol $p$ of $P$ in local coordinates near the boundary satisfies:
$$\partial^k_{x_n}\partial^\alpha_{\xi'}p_j(x',0,0,+1)=
(-1)^{j-|\alpha|}\partial^k_{x_n}\partial^\alpha_{\xi'}p_j(x',0,0,-1),$$
then $\pi^+P:C^{\infty}(X,{E_1})\rightarrow C^{\infty}(X,{E_2})$
by Section 2.1 of \cite{Wa1}.

Let $M$ be a compact manifold with boundary $\partial M$. We assume that the metric $g^{M}$ on $M$ has
the following form near the boundary
 \begin{equation}
 g^{M}=\frac{1}{h(x_{n})}g^{\partial M}+\texttt{d}x _{n}^{2} ,
\end{equation}
where $g^{\partial M}$ is the metric on $\partial M$. Let $U\subset
M$ be a collar neighborhood of $\partial M$ which is diffeomorphic $\partial M\times [0,1)$. By the definition of $h(x_n)\in C^{\infty}([0,1))$
and $h(x_n)>0$, there exists $\tilde{h}\in C^{\infty}((-\varepsilon,1))$ such that $\tilde{h}|_{[0,1)}=h$ and $\tilde{h}>0$ for some
sufficiently small $\varepsilon>0$. Then there exists a metric $\hat{g}$ on $\hat{M}=M\bigcup_{\partial M}\partial M\times
(-\varepsilon,0]$ which has the form on $U\bigcup_{\partial M}\partial M\times (-\varepsilon,0 ]$
 \begin{equation}
\hat{g}=\frac{1}{\tilde{h}(x_{n})}g^{\partial M}+\texttt{d}x _{n}^{2} ,
\end{equation}
such that $\hat{g}|_{M}=g$.
We fix a metric $\hat{g}$ on the $\hat{M}$ such that $\hat{g}|_{M}=g$.
Now we recall the main theorem in \cite{FGLS}.
\begin{thm}\label{th:32}{\bf(Fedosov-Golse-Leichtnam-Schrohe)}
 Let $X$ and $\partial X$ be connected, ${\rm dim}X=n\geq3$,
 $A=\left(\begin{array}{lcr}\pi^+P+G &   K \\
T &  S    \end{array}\right)$ $\in \mathcal{B}$ , and denote by $p$, $b$ and $s$ the local symbols of $P,G$ and $S$ respectively.
 Define:
 \begin{align}
{\rm{\widetilde{Wres}}}(A)=&\int_X\int_{\bf S}{\mathrm{Tr}}_E\left[p_{-n}(x,\xi)\right]\sigma(\xi)dx \nonumber\\
&+2\pi\int_ {\partial X}\int_{\bf S'}\left\{{\mathrm{Tr}}_E\left[({\mathrm{Tr}}b_{-n})(x',\xi')\right]+{\mathrm{Tr}}
_F\left[s_{1-n}(x',\xi')\right]\right\}\sigma(\xi')dx',
\end{align}
Then~~ a) ${\rm \widetilde{Wres}}([A,B])=0 $, for any
$A,B\in\mathcal{B}$;~~ b) It is a unique continuous trace on
$\mathcal{B}/\mathcal{B}^{-\infty}$.
\end{thm}

  Denote by $\sigma_{l}(H_{X})$ the $l$-order symbol of equivariant twisted Bismut Laplacian  $H_{X}$. An application of (2.1.4) in \cite{Wa1} shows that
\begin{equation}
\widetilde{{\rm Wres}}[\pi^+H_{X}^{-p_1}\circ\pi^+H_{X}^{-p_2}]=\int_M\int_{|\xi|=1}{\rm
Tr}_{S(TM)\otimes F}[\sigma_{-n}(H_{X}^{-p_1-p_2})]\sigma(\xi)\texttt{d}x+\int_{\partial
M}\Phi,
\end{equation}
where
 \begin{eqnarray}
\Phi&=&\int_{|\xi'|=1}\int^{+\infty}_{-\infty}\sum^{\infty}_{j, k=0}
\sum\frac{(-i)^{|\alpha|+j+k+1}}{\alpha!(j+k+1)!}
 {\rm Tr}_{S(TM)\otimes F}
\Big[\partial^j_{x_n}\partial^\alpha_{\xi'}\partial^k_{\xi_n}
\sigma^+_{r}(H_{X}^{-\frac{p_1}{2}})(x',0,\xi',\xi_n)\nonumber\\
&&\times\partial^\alpha_{x'}\partial^{j+1}_{\xi_n}\partial^k_{x_n}\sigma_{l}
(H_{X}^{-\frac{p_2}{2}})(x',0,\xi',\xi_n)\Big]d\xi_n\sigma(\xi')\texttt{d}x',
\end{eqnarray}
and the sum is taken over $r-k+|\alpha|+\ell-j-1=-n,r\leq-p_{1},\ell\leq-p_{2}$.

Let $M$ be a $n=\overline{n}+2$-dimensional compact oriented spin manifold with boundary
$\partial M$ and $n$ is an even integer. Since $[\sigma_{-n}(H_{X}^{\frac{-\overline{n}+2}{2}})]|_M$
has the same expression as $\sigma_{-n}(H_{X}^{\frac{-\overline{n}+2}{2}})$ in the case of
manifolds without boundary, then we only need to compute $\int_{\partial M}\Phi$.

Write
  \begin{equation}
D_x^{\alpha}=(-\sqrt{-1})^{|\alpha|}\partial_x^{\alpha};~\sigma(H_{X})=p_2+p_1+p_0;
~\sigma(H_{X}^{-1})=\sum^{\infty}_{j=2}q_{-j}.
\end{equation}
By the composition formula of psudodifferential operators, then we have
\begin{eqnarray}
1=\sigma(H_{X}\circ H_{X}^{-1})&=&\sum_{\alpha}\frac{1}{\alpha!}\partial^{\alpha}_{\xi}[\sigma(H_{X})]D^{\alpha}_{x}[\sigma(H_{X}^{-1})]\nonumber\\
&=&(p_2+p_1+p_0)(q_{-2}+q_{-3}+q_{-4}+\cdots)\nonumber\\
& &~~~+\sum_j(\partial_{\xi_j}p_2+\partial_{\xi_j}p_1+\partial_{\xi_j}p_0)(
D_{x_j}q_{-2}+D_{x_j}q_{-3}+D_{x_j}q_{-4}+\cdots) .
\end{eqnarray}
Then we have
 \begin{equation}
\sigma_{2}(H_{X}^{-1})=|\xi|^{2}~~ \sigma_{-2}(H_{X}^{-1})=|\xi|^{-2},
~~\sigma_{-(\overline{n}-2)}(H_{X}^{\frac{-\overline{n}+2}{2}})=(|\xi|^2)^{1-\frac{\overline{n}}{2}}.
\end{equation}
 \begin{align}\label{eq:3.12}
\sigma_{-3}(H_{X}^{-1})=&-\sqrt{-1}|\xi|^{-4}\xi_k(\Gamma^k-2\delta^k)-\sqrt{-1}|\xi|^{-6}2\xi^j\xi_\alpha\xi_\beta
\partial_jg^{\alpha\beta}-\sum_{j=1}^{n}\sqrt{-1}|\xi|^{-4}(X_{i}-\frac{1}{2} \langle X,\partial_{j}  \rangle)\xi_{j}\nonumber\\
&-\sum_{j=1}^{n}\sqrt{-1}|\xi|^{-4}\Big( c(\partial_{j})\big((c(T)+c(Y)) \otimes \Phi
 +\big(c(T)+c(Y)\big)\otimes \Phi c(\partial_{j})\Big)\xi_{j}; \\
\sigma_{-(\overline{n}-2)}(H_{X}^{\frac{-\overline{n}+2}{2}})&=(|\xi|^2)^{1-\frac{\overline{n}}{2}}.
\end{align}
Then we  obtain by induction
\begin{eqnarray}
\sigma_{1-\overline{n}}(H_{X}^{\frac{-\overline{n}+2}{2}})(x,\xi)
            &=& \sigma_{3-\overline{n}}(H_{X}^{\frac{-\overline{n}+4}{2}})\sigma_{2}^{-1}
            +\sigma_{2}^{(-\frac{\overline{n}}{2}+2)}\sigma_{-3}(H_{X}^{-1})
                  -\sqrt{-1}\partial_{\xi_{\mu}}\sigma_{2}^{(-\frac{\overline{n}}{2}+2)}\partial_{x_{\mu}}\sigma_{2}^{-1}\nonumber\\
          &=& \frac{\overline{n}-2}{2}\sigma_{2}^{(-\frac{\overline{n}}{2}+2)}\sigma_{-3}(H_{X}^{-1})-\sqrt{-1} \sum_{k=0}^{\frac{\overline{n}}{2}-3}
             \partial_{\xi_{\mu}}\sigma_{2}^{(-\frac{\overline{n}}{2}+k+2)}\partial_{x_{\mu}}\sigma_{2}^{-1}\big(\sigma_{2}^{-1}\big)^{k}.
\end{eqnarray}

 \begin{lem}\label{le:31}
The symbols of the equivariant twisted Bismut Laplacian with torsion $H_{X}$:
 \begin{align}
\sigma_{-2}(H_{X}^{-1})=&|\xi|^{-2}; \\
\sigma_{-(\overline{n}-2)}(H_{X}^{\frac{-\overline{n}+2}{2}})&=(|\xi|^2)^{1-\frac{\overline{n}}{2}}; \\
\sigma_{1-\overline{n}}(D^{-\overline{n}+2})=&\frac{\overline{n}-2}{2}\sigma_{2}^{(-\frac{\overline{n}}{2}+2)}\sigma_{-3}(D^{-2})-\sqrt{-1} \sum_{k=0}^{\frac{\overline{n}}{2}-3}
             \partial_{\xi_{\mu}}\sigma_{2}^{(-\frac{\overline{n}}{2}+k+2)}\partial_{x_{\mu}}\sigma_{2}^{-1}\big(\sigma_{2}^{-1}\big)^{k};\\
\sigma_{1-\overline{n}}(H_{X}^{\frac{-\overline{n}+2}{2}})=&\sigma_{1-\overline{n}}(D^{-\overline{n}+2})
-\frac{\overline{n}-2}{2}\sigma_{2}^{(-\frac{\overline{n}}{2}+2)}\sum_{j=1}^{n}\sqrt{-1}|\xi|^{-4}(X_{i}
-\frac{1}{2} \langle X,\partial_{j}  \rangle)\xi_{j}\nonumber\\
&-\frac{\overline{n}-2}{2}\sigma_{2}^{(-\frac{\overline{n}}{2}+2)}\sum_{j=1}^{n}\sqrt{-1}|\xi|^{-4}\Big( c(\partial_{j})\big((c(T)+c(Y)) \otimes \Phi
 +\big(c(T)+c(Y)\big)\otimes \Phi c(\partial_{j})\Big)\xi_{j}.
\end{align}
\end{lem}

Since $\Phi$ is a global form on $\partial M$, so for any fixed point $x_{0}\in\partial M$, we can choose the normal coordinates
$U$ of $x_{0}$ in $\partial M$(not in $M$) and compute $\Phi(x_{0})$ in the coordinates $\widetilde{U}=U\times [0,1)$ and the metric
$\frac{1}{h(x_{n})}g^{\partial M}+\texttt{d}x _{n}^{2}$. The dual metric of $g^{M}$ on $\widetilde{U}$ is
$h(x_{n})g^{\partial M}+\texttt{d}x _{n}^{2}.$ Write
$g_{ij}^{M}=g^{M}(\frac{\partial}{\partial x_{i}},\frac{\partial}{\partial x_{j}})$;
$g^{ij}_{M}=g^{M}(d x_{i},dx_{j})$, then

\begin{equation*}
[g_{i,j}^{M}]=
\begin{bmatrix}\frac{1}{h( x_{n})}[g_{i,j}^{\partial M}]&0\\0&1\end{bmatrix};\quad
[g^{i,j}_{M}]=\begin{bmatrix} h( x_{n})[g^{i,j}_{\partial M}]&0\\0&1\end{bmatrix},
\end{equation*}
and
\begin{equation*}
\partial_{x_{s}} g_{ij}^{\partial M}(x_{0})=0,\quad 1\leq i,j\leq n-1;\quad g_{i,j}^{M}(x_{0})=\delta_{ij}.
\end{equation*}

Let $\{E_{1},\cdots, E_{n-1}\}$ be an orthonormal frame field in $U$ about $g^{\partial M}$ which is parallel along geodesics and
$E_{i}=\frac{\partial}{\partial x_{i}}(x_{0})$, then $\{\widetilde{E_{1}}=\sqrt{h(x_{n})}E_{1}, \cdots,
\widetilde{E_{n-1}}=\sqrt{h(x_{n})}E_{n-1},\widetilde{E_{n}}=dx_{n}\}$ is the orthonormal frame field in $\widetilde{U}$ about $g^{M}.$
Locally $S(TM)|\widetilde{U}\cong \widetilde{U}\times\wedge^{*}_{C}(\frac{n}{2}).$ Let $\{f_{1},\cdots,f_{n}\}$ be the orthonormal basis of
$\wedge^{*}_{C}(\frac{n}{2})$. Take a spin frame field $\sigma: \widetilde{U}\rightarrow Spin(M)$ such that
$\pi\sigma=\{\widetilde{E_{1}},\cdots, \widetilde{E_{n}}\}$ where $\pi: Spin(M)\rightarrow O(M)$ is a double covering, then
$\{[\sigma, f_{i}], 1\leq i\leq n\}$ is an orthonormal frame of $S(TM)|_{\widetilde{U}}.$ In the following, since the global form $\Phi$
is independent of the choice of the local frame, so we can compute $\texttt{tr}_{S(TM)}$ in the frame $\{[\sigma, f_{i}], 1\leq i\leq n\}$.
Let $\{\hat{E}_{1},\cdots,\hat{E}_{n}\}$ be the canonical basis of $R^{n}$ and
$c(\hat{E}_{i})\in cl_{C}(n)\cong Hom(\wedge^{*}_{C}(\frac{n}{2}),\wedge^{*}_{C}(\frac{n}{2}))$ be the Clifford action. Then
\begin{equation*}
c(\widetilde{E_{i}})=[(\sigma,c(\hat{E}_{i}))]; \quad c(\widetilde{E_{i}})[(\sigma, f_{i})]=[\sigma,(c(\hat{E}_{i}))f_{i}]; \quad
\frac{\partial}{\partial x_{i}}=[(\sigma,\frac{\partial}{\partial x_{i}})],
\end{equation*}
then we have $\frac{\partial}{\partial x_{i}}c(\widetilde{E_{i}})=0$ in the above frame. By Lemma 2.2 in \cite{Wa3}, we have

\begin{lem}\label{le:32}
With the metric $g^{M}$ on $M$ near the boundary
\begin{eqnarray}
\partial_{x_j}(|\xi|_{g^M}^2)(x_0)&=&\left\{
       \begin{array}{c}
        0,  ~~~~~~~~~~ ~~~~~~~~~~ ~~~~~~~~~~~~~{\rm if }~j<n; \\[2pt]
       h'(0)|\xi'|^{2}_{g^{\partial M}},~~~~~~~~~~~~~~~~~~~~~{\rm if }~j=n.
       \end{array}
    \right. \\
\partial_{x_j}[c(\xi)](x_0)&=&\left\{
       \begin{array}{c}
      0,  ~~~~~~~~~~ ~~~~~~~~~~ ~~~~~~~~~~~~~{\rm if }~j<n;\\[2pt]
\partial x_{n}(c(\xi'))(x_{0}), ~~~~~~~~~~~~~~~~~{\rm if }~j=n,
       \end{array}
    \right.
\end{eqnarray}
where $\xi=\xi'+\xi_{n}\texttt{d}x_{n}$.
\end{lem}

From the remark above, now we can compute $\Phi$ (see formula (3.48) for the definition of $\Phi$).
Since the sum is taken over $r+\ell-k-j-|\alpha|-1=-(\overline{n}+2), \ r\leq-2, \ell\leq 2-\overline{n}$, then we have the $\int_{\partial_{M}}\Phi$
is the sum of the following five cases:

\noindent  {\bf case a)~I)}~$r=-2,~l=2-\overline{n},~k=j=0,~|\alpha|=1$

 By  (3.48), we get
 \begin{equation}
{\rm case~a)~I)}=-\int_{|\xi'|=1}\int^{+\infty}_{-\infty}\sum_{|\alpha|=1}{\rm{Tr}}
\Big[\partial^\alpha_{\xi'}\pi^+_{\xi_n}\sigma_{-2}(H_{X}^{-1})\times\partial^\alpha_{x'}
\partial_{\xi_n}\sigma_{-(\overline{n}-2)}(H_{X}^{\frac{-(\overline{n}-2)}{2}})\Big](x_0)d\xi_n\sigma(\xi')dx'.
\end{equation}
By Lemma 3.6 and Lemma 3.7, for $i<n$, then
 \begin{equation}
\partial_{x_i}\sigma_{-(\overline{n}-2)}(H_{X}^{\frac{-(\overline{n}-2)}{2}})(x_0)=\partial_{x_i}{(|\xi|^{(2-\overline{n})})}(x_0)
 =(1-\frac{\overline{n}}{2}) (|\xi|^2)^{-\frac{\overline{n}}{2}}
 \partial_{x_i}(|\xi|^2)(x_0)=0,
\end{equation}
 so {\rm {\bf case~a)~I)}} vanishes.

\noindent  {\bf case a)~II)}~$r=-2,~l=2-\overline{n},~k=|\alpha|=0,~j=1$

 By (3.48), we get
  \begin{equation}
{\rm case~a)~II)}=-\frac{1}{2}\int_{|\xi'|=1}\int^{+\infty}_{-\infty}{\rm{Tr}}
 \Big[\partial_{x_n}\pi^+_{\xi_n}\sigma_{-2}(H_{X}^{-1})\times \partial_{\xi_n}^2
\sigma_{-(\overline{n}-2)}(H_{X}^{\frac{-(\overline{n}-2)}{2}})\Big](x_0)d\xi_n\sigma(\xi')dx'.
\end{equation}
By Lemma 3.6, we have
\begin{equation}
\partial_{x_n}\sigma_{-2}(H_{X}^{-1})(x_0)|_{|\xi'|=1}=-\frac{h'(0)}{(1+\xi_n^2)^2}.
\end{equation}
By the Cauchy integral formula, then
\begin{eqnarray}
\pi^+_{\xi_n}\partial_{x_n}\sigma_{-2}(H_{X}^{-1})(x_0)|_{|\xi'|=1}
&=&-h'(0)\frac{1}{2\pi i} \lim_{u\rightarrow 0^-}\int_{\Gamma^+}\frac{\frac{1}{(\eta_n+i)^2(\xi_n+iu-\eta_n)}}{(\eta_n-i)^2}d\eta_n\nonumber\\
&=&h'(0)\frac{i\xi_n+2}{4(\xi_n-i)^2}.
\end{eqnarray}
From  Lemma 3.6  we obtain
\begin{eqnarray}
\partial^2_{\xi_n} \big(\sigma_{-(\overline{n}-2)}(H_{X}^{\frac{-(\overline{n}-2)}{2}})\big)(x_0)
&=&\partial^2_{\xi_n}\big((|\xi|^2)^{1-\frac{\overline{n}}{2}}\big) (x_0)
\nonumber\\
&=& (1-\frac{\overline{n}}{2})(-\frac{\overline{n}}{2})(|\xi|^{2})^{-\frac{\overline{n}}{2}-1} \big(\partial_{\xi_n}|\xi|^2\big)^{2}(x_0)\nonumber\\
&&+(1-\frac{\overline{n}}{2})(|\xi|^{2})^{-\frac{\overline{n}}{2}} \partial^2_{\xi_n}(|\xi|^2 (x_0))\nonumber\\
&=& \Big((2\overline{n}-2) \xi_n^{2}-2 \Big)(\frac{\overline{n}}{2}-1)(1+\xi_n^{2})^{(-\frac{\overline{n}}{2}-1)}.
\end{eqnarray}
Then
\begin{eqnarray}
&& \int_{-\infty}^{\infty}h'(0)\frac{i\xi_n+2}{4(\xi_n-i)^2}\times
 \Big((2\overline{n}-2) \xi_n^{2}-2 \Big)(\frac{\overline{n}}{2}-1)(1+\xi_n^{2})^{(-\frac{\overline{n}}{2}-1)}d\xi_n\nonumber\\
&=&\frac{h'(0)}{4}(\frac{\overline{n}}{2}-1)   \int_{\Gamma^+}\frac{(2\overline{n}-2)i\xi_n^3+(4\overline{n}-4)\xi_n^2-2i\xi_n-4}
   {(\xi_n-i)^{(\frac{\overline{n}}{2}+3)}(\xi_n+i)^{(\frac{\overline{n}}{2}+1)}}d\xi_n\nonumber\\
&=& \frac{h'(0)}{4}(\frac{\overline{n}}{2}-1)  \frac{2\pi i}{(\frac{\overline{n}}{2}+2)!}
\left[\frac{(2\overline{n}-2)i\xi_n^3+(4\overline{n}-4)\xi_n^2-2i\xi_n-4}{(\xi_n+i)^{(\frac{\overline{n}}{2}+1)}}\right]
 ^{(\frac{\overline{n}}{2}+2)}|_{\xi_n=i}\nonumber\\
 &:=& \frac{h'(0)}{4} (\frac{\overline{n}}{2}-1) \frac{2\pi i}{(\frac{\overline{n}}{2}+2)!}  L_{0}.
\end{eqnarray}
Since $n=\overline{n}+2$ is even, ${\rm tr}_{S(TM)}[{\rm id}]={\rm dim}(\wedge^*(\frac{\overline{n}+2}{2} ))=2^{(\frac{\overline{n}}{2}+1)}.$
Then we obtain
 \begin{equation}
{\rm {\bf case~a)~II)}}=Vol(S_{\overline{n}}) \frac{-h'(0)}{4} (\frac{\overline{n}}{2}-1)
\frac{\pi i}{(\frac{\overline{n}}{2}+2)!} 2^{(\frac{\overline{n}}{2}+1)} L_{0}dx',
\end{equation}
 where $Vol(S_{\overline{n}})$ is the canonical volume of $S_{\overline{n}}$

\noindent  {\bf case a)~III)}~$r=-2,~l=2-\overline{n},~j=|\alpha|=0,~k=1$

 By (3.48) and an integration by parts, we get
 \begin{eqnarray}
{\rm case~ a)~III)}&=&-\frac{1}{2}\int_{|\xi'|=1}\int^{+\infty}_{-\infty}{\rm trace} \Big[\partial_{\xi_n}\pi^+_{\xi_n}\sigma_{-2}(H_{X}^{-1})\times
                      \partial_{\xi_n}\partial_{x_n}\sigma_{-(\overline{n}-2)}(H_{X}^{\frac{-(\overline{n}-2)}{2}})\Big](x_0)d\xi_n\sigma(\xi')dx'\nonumber\\
                   &=& \frac{1}{2}\int_{|\xi'|=1}\int^{+\infty}_{-\infty} {\rm trace}\Big[\partial_{\xi_n}^2\pi^+_{\xi_n}\sigma_{-2}(H_{X}^{-1})\times
                       \partial_{x_n}\sigma_{-(\overline{n}-2)}(H_{X}^{\frac{-(\overline{n}-2)}{2}})\Big](x_0)d\xi_n\sigma(\xi')dx'.
\end{eqnarray}
 By Lemma 3.6, we have
\begin{equation}
\partial_{\xi_n}^2\pi_{\xi_n}^+\sigma_{-2}(H_{X}^{-1})(x_0)|_{|\xi'|=1}=\frac{-i}{(\xi_n-i)^3}.
\end{equation}
And
\begin{equation}
\partial_{x_n} \big(\sigma_{-(\overline{n}-2)}(H_{X}^{\frac{-(\overline{n}-2}{2})})\big)(x_0)
=\partial_{x_n}\big((|\xi|^2)^{1-\frac{\overline{n}}{2}}\big) (x_0)
=h'(0)(1-\frac{\overline{n}}{2})(1+\xi_{n}^{2})^{-\frac{\overline{n}}{2}}.
\end{equation}
Then
\begin{eqnarray}
&& \int_{-\infty}^{\infty}{\rm{Tr}}\Big[\frac{-i}{(\xi_n-i)^3}\times h'(0)(1-\frac{\overline{n}}{2})(1+\xi_{n}^{2})^{-\frac{\overline{n}}{2}}
\Big] d\xi_n\nonumber\\
&=&-i h'(0)(1-\frac{\overline{n}}{2}) 2^{(\frac{\overline{n}}{2}+1)}
      \int_{\Gamma^+} \frac{1}
   {(\xi_n-i)^{(\frac{\overline{n}}{2}+3)}(\xi_n+i)^{(\frac{\overline{n}}{2} )}}d\xi_n\nonumber\\
&=& -i h'(0)(1-\frac{\overline{n}}{2}) 2^{(\frac{\overline{n}}{2}+1)} \frac{2\pi i}{(\frac{\overline{n}}{2}+2)!}
\left[\frac{1}{(\xi_n+i)^{(\frac{\overline{n}}{2})}}\right]
 ^{(\frac{\overline{n}}{2}+2)}|_{\xi_n=i}\nonumber\\
 &=& -i h'(0)(1-\frac{\overline{n}}{2}) 2^{(\frac{\overline{n}}{2}+1)} \frac{2\pi i}{(\frac{\overline{n}}{2}+2)!}  L_{1}.
\end{eqnarray}
Therefore
  \begin{equation}
{\rm {\bf case~a)~III)}}=-iVol(S_{\overline{n}})h'(0)(1-\frac{\overline{n}}{2}) 2^{(\frac{\overline{n}}{2}+1)}
\frac{\pi i}{(\frac{\overline{n}}{2}+2)!}  L_{1}dx'.
\end{equation}

 \noindent  {\bf case b)}~$r=-2,~l=1-\overline{n},~k=j=|\alpha|=0$

 By (3.48) and an integration by parts, we get
 \begin{eqnarray}
{\rm case~ b)}&=&-i\int_{|\xi'|=1}\int^{+\infty}_{-\infty}{\rm trace} [\pi^+_{\xi_n}\sigma_{-2}(H_{X}^{-1})\times
     \partial_{\xi_n}\sigma_{1-\overline{n}}(H_{X}^{\frac{-\overline{n}+2})}{2}](x_0)d\xi_n\sigma(\xi')dx'\nonumber\\
     &=&i\int_{|\xi'|=1}\int^{+\infty}_{-\infty}{\rm trace} [\partial_{\xi_n}\pi^+_{\xi_n}\sigma_{-2}(H_{X}^{-1})\times
        \sigma_{1-\overline{n}}(H_{X}^{\frac{-\overline{n}+2}{2}})](x_0)d\xi_n\sigma(\xi')dx'.
\end{eqnarray}
 By Lemma 3.6, we have
\begin{equation}
\partial_{\xi_n}\pi_{\xi_n}^+\sigma_{-2}(H_{X}^{-1})(x_0)|_{|\xi'|=1}=\frac{i}{2(\xi_n-i)^2}.
\end{equation}
Let
\begin{equation}
\sigma_{1-\overline{n}}(H_{X}^{\frac{-\overline{n}+2}{2}})£º=A_{1}+A_{2}+A_{3},
\end{equation}
where
 \begin{align}
A_{1}=&\sigma_{1-\overline{n}}(D^{-\overline{n}+2});\\
A_{2}=& -\frac{\overline{n}-2}{2}\sigma_{2}^{(-\frac{\overline{n}}{2}+2)}\sum_{j=1}^{n}\sqrt{-1}|\xi|^{-4}(X_{i}
-\frac{1}{2} \langle X,\partial_{j}  \rangle)\xi_{j};\\
A_{3}=&-\frac{\overline{n}-2}{2}\sigma_{2}^{(-\frac{\overline{n}}{2}+2)}\sum_{j=1}^{n}\sqrt{-1}|\xi|^{-4}\Big( c(\partial_{j})\big((c(T)+c(Y)) \otimes \Phi
 +\big(c(T)+c(Y)\big)\otimes \Phi c(\partial_{j})\Big)\xi_{j}.
\end{align}
From (3.75) and (3.77), we obtain
\begin{eqnarray}
&& i\int_{|\xi'|=1}\int^{+\infty}_{-\infty}{\rm trace}
\Big\{\frac{i}{2(\xi_n-i)^2}\times   A_{1} \Big\}d\xi_n\sigma(\xi')dx' \nonumber\\
&=& i\int_{|\xi'|=1}\int^{+\infty}_{-\infty}{\rm trace}
\Bigg\{\frac{i}{2(\xi_n-i)^2}\times \Bigg[\frac{\overline{n}-2}{2} (1+\xi_{n}^{2})^{(-\frac{\overline{n}}{2}+2)}
   \Big(\frac{-i}{(1+\xi_n^2)^2}\times\frac{\overline{n}+1}{2}h'(0)\xi_n\nonumber\\
   &&-\frac{2ih'(0)\xi_n}{(1+\xi_n^2)^3} \Big)+\sqrt{-1} h'(0)(-\frac{\overline{n}^{2}}{4}
   +\frac{3\overline{n}}{2}-2)  \xi_{n} (1+\xi_{n}^{2})^{(-\frac{\overline{n}}{2}-1)}\Bigg]
\Bigg\}d\xi_n\sigma(\xi')dx' \nonumber\\
&=&-\frac{i h'(0)}{8} Vol(S_{\overline{n}})\int_{\Gamma^+}\frac{-(\overline{n}-2)(\overline{n}+1)\xi_{n}^{3}+(-2\overline{n}^{2}
   +3\overline{n}-2)  \xi_{n}}{(\xi_n-i)^{(\frac{\overline{n}}{2}+3)}(\xi_n+i)^{(\frac{\overline{n}}{2}+1)}}d\xi_ndx'\nonumber\\
&=&-\frac{i h'(0)}{8} Vol(S_{\overline{n}})2^{(\frac{\overline{n}}{2}+1)}  \frac{2\pi i}{(\frac{\overline{n}}{2}+2)!}
  \left[\frac{-(\overline{n}-2)(\overline{n}+1)\xi_{n}^{3}+(-2\overline{n}^{2}
   +3\overline{n}-2)  \xi_{n}}{(\xi_n+i)^{(\frac{\overline{n}}{2}+1)}}\right]^{(\frac{\overline{n}}{2}+2)}|_{\xi_n=i}dx'\nonumber\\
&=&  \frac{\pi h'(0) Vol(S_{\overline{n}})2^{(\frac{\overline{n}}{2}-1)}}{(\frac{\overline{n}}{2}+2)!}
  \left[\frac{-(\overline{n}-2)(\overline{n}+1)\xi_{n}^{3}+(-2\overline{n}^{2}
   +3\overline{n}-2)  \xi_{n}}{(\xi_n+i)^{(\frac{\overline{n}}{2}+1)}}\right]^{(\frac{\overline{n}}{2}+2)}|_{\xi_n=i}dx'\nonumber\\
&=&  \frac{\pi h'(0) Vol(S_{\overline{n}})2^{(\frac{\overline{n}}{2}-1)}}{(\frac{\overline{n}}{2}+2)!}L_{2}dx'.
\end{eqnarray}
From (3.75) and (3.78), we obtain
\begin{eqnarray}
&& i\int_{|\xi'|=1}\int^{+\infty}_{-\infty}{\rm trace}
\Big\{\frac{i}{2(\xi_n-i)^2}\times   A_{2} \Big\}d\xi_n\sigma(\xi')dx' \nonumber\\
&=& i\int_{|\xi'|=1}\int^{+\infty}_{-\infty}{\rm trace}
\Bigg\{\frac{i}{2(\xi_n-i)^2}\times \Bigg[-\frac{\overline{n}-2}{2}|\xi|^{- \overline{n}+4}\sum_{j=1}^{n}\sqrt{-1}|\xi|^{-4}(X_{i}
-\frac{1}{2} \langle X,\partial_{j}  \rangle)\xi_{j}\Bigg]
\Bigg\}d\xi_n\sigma(\xi')dx' \nonumber\\
&=&  \frac{(2-\overline{n})2^{(\frac{\overline{n}}{2}-2)}\pi X_{n} Vol(S_{\overline{n}})}{(\frac{\overline{n}}{2}+1)!}
  \Big[\frac{\xi_n }{(\xi_n+i)^{(\frac{\overline{n}}{2})}}\Big]^{(\frac{\overline{n}}{2}+1)}|_{\xi_n=i}dx'.
\end{eqnarray}
Similarly,
\begin{eqnarray}
&& i\int_{|\xi'|=1}\int^{+\infty}_{-\infty}{\rm trace}
\Big\{\frac{i}{2(\xi_n-i)^2}\times   A_{3} \Big\}d\xi_n\sigma(\xi')dx' \nonumber\\
&=& i\int_{|\xi'|=1}\int^{+\infty}_{-\infty}{\rm trace}
\Bigg\{-\frac{\overline{n}-2}{2}\sigma_{2}^{(-\frac{\overline{n}}{2}+2)}\sum_{j=1}^{n}\sqrt{-1}|\xi|^{-4}\Big( c(\partial_{j})\big((c(T)+c(Y))
 \otimes \Phi\nonumber\\
&& +\big(c(T)+c(Y)\big)\otimes \Phi c(\partial_{j})\Big)\xi_{j}
\Bigg\}d\xi_n\sigma(\xi')dx' \nonumber\\
&=&  \frac{(2-\overline{n})2^{(\frac{\overline{n}}{2}-2)}\pi (-2g(Y, \partial_{n})) Vol(S_{\overline{n}})}{(\frac{\overline{n}}{2}+1)!}
  \Big[\frac{\xi_n }{(\xi_n+i)^{(\frac{\overline{n}}{2})}}\Big]^{(\frac{\overline{n}}{2}+1)}|_{\xi_n=i}dx'.
\end{eqnarray}
Therefore
\begin{eqnarray}
{\bf case~ b)}&=& \frac{\pi h'(0) Vol(S_{\overline{n}})2^{(\frac{\overline{n}}{2}-1)}}{(\frac{\overline{n}}{2}+2)!}L_{2}dx'\nonumber\\
&&+ \frac{(2-\overline{n})2^{(\frac{\overline{n}}{2}-2)}\pi(X_{n}-2g(Y, \partial_{n})) Vol(S_{\overline{n}})}{(\frac{\overline{n}}{2}+1)!}
  \Big[\frac{\xi_n }{(\xi_n+i)^{(\frac{\overline{n}}{2})}}\Big]^{(\frac{\overline{n}}{2}+1)}|_{\xi_n=i}dx'
\end{eqnarray}

\noindent {\bf  case c)}~$r=-3,~l=2-\overline{n},~k=j=|\alpha|=0$

 By (3.48), we get
\begin{equation}
{\rm case~ c)}=-i\int_{|\xi'|=1}\int^{+\infty}_{-\infty}{\rm trace} \Big[\pi^+_{\xi_n}\sigma_{-3}(H_{X}^{-1})\times
\partial_{\xi_n}\sigma_{-(\overline{n}-2)}(H_{X}^{\frac{-(\overline{n}-2)}{2}})\Big](x_0)d\xi_n\sigma(\xi')dx'.
\end{equation}
From Lemma 3.6, we get
\begin{equation}
\partial_{\xi_n}\sigma_{-(\overline{n}-2)}(H_{X}^{\frac{-(\overline{n}-2)}{2}})(x_0)
=\partial_{\xi_n}((|\xi|^2)^{1-\frac{\overline{n}}{2}})(x_0)=2(1-\frac{\overline{n}}{2})\xi_n(1+\xi_n^2)^{-\frac{\overline{n}}{2}}.
\end{equation}
Let
 \begin{align}\label{eq:3.12}
 \sigma_{-3}(H_{X}^{-1})=B_{1}+B_{2}+B_{3},
\end{align}
where
 \begin{align}\label{eq:3.12}
B_{1}=&-\sqrt{-1}|\xi|^{-4}\xi_k(\Gamma^k-2\delta^k)-\sqrt{-1}|\xi|^{-6}2\xi^j\xi_\alpha\xi_\beta
\partial_jg^{\alpha\beta};  \\
 B_{2}=& -\sum_{j=1}^{n}\sqrt{-1}|\xi|^{-4}(X_{i}-\frac{1}{2} \langle X,\partial_{j}  \rangle)\xi_{j};  \\
 B_{3}=&-\sum_{j=1}^{n}\sqrt{-1}|\xi|^{-4}\Big( c(\partial_{j})\big((c(T)+c(Y)) \otimes \Phi
 +\big(c(T)+c(Y)\big)\otimes \Phi c(\partial_{j})\Big)\xi_{j}.
\end{align}
From (3.85) and (3.87), we obtain
\begin{eqnarray}
&&-i\int_{|\xi'|=1}\int^{+\infty}_{-\infty}{\rm trace} \Big[\pi^+_{\xi_n}(B_{1})\times
\partial_{\xi_n}\sigma_{-(\overline{n}-2)}(H_{X}^{\frac{-(\overline{n}-2)}{2}})\Big](x_0)d\xi_n\sigma(\xi')dx' \nonumber\\
&=&-i\int_{|\xi'|=1}\int^{+\infty}_{-\infty}{\rm trace}
\Bigg[ \frac{2(1-\frac{\overline{n}}{2})\xi_n}{(1+\xi_n^2)^{\frac{\overline{n}}{2}} }
\times i h'(0)\Big(\frac{i\overline{n}}{8(\xi_n-i)^2}+ \frac{1}{4(\xi_n-i)^3} \Big)
\Bigg]d\xi_n\sigma(\xi')dx' \nonumber\\
&=& (1-\frac{\overline{n}}{2})2^{(\frac{\overline{n}}{2}-1)} Vol(S_{\overline{n}})h'(0)
  \int_{\Gamma^+}\frac{(i\overline{n}\xi_{n} + \overline{n}+2)\xi_{n}}
  {(\xi_n+i)^{(\frac{\overline{n}}{2})}(\xi_n-i)^{(\frac{\overline{n}}{2}+3)}}d\xi_ndx' \nonumber\\
&=& (1-\frac{\overline{n}}{2})2^{(\frac{\overline{n}}{2}-1)} Vol(S_{\overline{n}})h'(0) \frac{2\pi i}{(\frac{\overline{n}}{2}+2)!}
  \left[\frac{(i\overline{n}\xi_{n} + \overline{n}+2)\xi_{n}}
 {(\xi_n+i)^{\frac{\overline{n}}{2}}}\right]^{(\frac{\overline{n}}{2}+2)}|_{\xi_n=i}dx'\nonumber\\
 & =& (1-\frac{\overline{n}}{2})2^{(\frac{\overline{n}}{2}-1)} Vol(S_{\overline{n}})h'(0) \frac{2\pi i}{(\frac{\overline{n}}{2}+2)!}
    L_{3}dx'.
\end{eqnarray}
Similarly,
\begin{eqnarray}
&&-i\int_{|\xi'|=1}\int^{+\infty}_{-\infty}{\rm trace} \Big[\pi^+_{\xi_n}(B_{2})\times
\partial_{\xi_n}\sigma_{-(\overline{n}-2)}(H_{X}^{\frac{-(\overline{n}-2)}{2}})\Big](x_0)d\xi_n\sigma(\xi')dx' \nonumber\\
&=&  \frac{(2-\overline{n})2^{(\frac{\overline{n}}{2}-2)}\pi X_{n} Vol(S_{\overline{n}})}{(\frac{\overline{n}}{2}+1)!}
  \Big[\frac{\xi_n }{(\xi_n+i)^{(\frac{\overline{n}}{2})}}\Big]^{(\frac{\overline{n}}{2}+1)}|_{\xi_n=i}dx',
\end{eqnarray}
and
\begin{eqnarray}
&&-i\int_{|\xi'|=1}\int^{+\infty}_{-\infty}{\rm trace} \Big[\pi^+_{\xi_n}(B_{3})\times
\partial_{\xi_n}\sigma_{-(\overline{n}-2)}(H_{X}^{\frac{-(\overline{n}-2)}{2}})\Big](x_0)d\xi_n\sigma(\xi')dx' \nonumber\\
&=&  \frac{(2-\overline{n})2^{(\frac{\overline{n}}{2}-2)}\pi (-2g(Y, \partial_{n})) Vol(S_{\overline{n}})}{(\frac{\overline{n}}{2}+1)!}
  \Big[\frac{\xi_n }{(\xi_n+i)^{(\frac{\overline{n}}{2})}}\Big]^{(\frac{\overline{n}}{2}+1)}|_{\xi_n=i}dx'.
\end{eqnarray}
Then
\begin{eqnarray}
{\bf case~ c)}&=&(1-\frac{\overline{n}}{2})2^{(\frac{\overline{n}}{2}-1)} Vol(S_{\overline{n}})h'(0) \frac{2\pi i}{(\frac{\overline{n}}{2}+2)!}
    L_{3}dx'\nonumber\\
    &&+ \frac{(2-\overline{n})2^{(\frac{\overline{n}}{2}-2)}\pi(X_{n}-2g(Y, \partial_{n}))  Vol(S_{\overline{n}})}{(\frac{\overline{n}}{2}+1)!}
  \Big[\frac{\xi_n }{(\xi_n+i)^{(\frac{\overline{n}}{2})}}\Big]^{(\frac{\overline{n}}{2}+1)}|_{\xi_n=i}dx'.
\end{eqnarray}

Now $\Phi$ is the sum of the cases a), b) and c), where
\begin{align}\label{eq:3.19}
  \Big[\frac{\xi_n }{(\xi_n+i)^{n}}\Big]^{(n+1)}|_{\xi_n=i}
  =& \frac{-(n+2)(n+3)(n+4)\cdots (2n+1)}{2^{2n+3}},   \\
  \Big[\frac{\xi_n^{2} }{(\xi_n+i)^{n}}\Big]^{(n+1)}|_{\xi_n=i}
  =& \frac{-(n-1)n(n+1)(n+2)\cdots (2n-1)}{2^{2n}i},  \\
    \Big[\frac{\xi_n ^{3}}{(\xi_n+i)^{n}}\Big]^{(n+1)}|_{\xi_n=i}
  =& \frac{3n(n+1)(n+2)\cdots (2n-1)}{2^{2n+1}},
\end{align}
then
\begin{eqnarray}
\Phi&=&  Vol(S_{\overline{n}})h'(0)(\frac{\overline{n}}{2}-1)\frac{2\pi i}{(\frac{\overline{n}}{2}+2)!} 2^{(\frac{\overline{n}}{2}-2)}
\Big[ \frac{   (-2i\overline{n}+4i)\xi_{n}^{3}-4\overline{n}\xi_{n}^{2}+(2i\overline{n}+4i)\xi_{n}}
 {(\xi_n+i)^{(\frac{\overline{n}}{2}+1)}}\Big]^{(\frac{\overline{n}}{2}+2)}\Big|_{\xi_n=i} \nonumber\\
    &&+ \frac{(3\overline{n}-6)2^{(\frac{\overline{n}}{2}-2)}\pi(X_{n}-2g(Y, \partial_{n}))  Vol(S_{\overline{n}})}{(\frac{\overline{n}}{2}+1)!}
  \Big[\frac{\xi_n }{(\xi_n+i)^{(\frac{\overline{n}}{2})}}\Big]^{(\frac{\overline{n}}{2}+1)}|_{\xi_n=i}dx'\nonumber\\
&=&  Vol(S_{\overline{n}})h'(0)(\frac{\overline{n}}{2}-1)\frac{2\pi i}{(\frac{\overline{n}}{2}+2)!}
\frac{(-\frac{1}{2}i\overline{n}^{4}-i \overline{n}^{3}+\frac{5}{4}i\overline{n}^{2}+\frac{5}{2}i\overline{n} )
(\frac{\overline{n}}{2}+2 )(\frac{\overline{n}}{2}+3 )\cdots (\overline{n}-1)
}{2^{\frac{\overline{n}}{2}+2}}
 \nonumber\\
    &&+ \frac{(3\overline{n}-6)2^{(\frac{\overline{n}}{2}-2)}\pi(X_{n}-2g(Y, \partial_{n}))  Vol(S_{\overline{n}})}{(\frac{\overline{n}}{2}+1)!}
  \Big(\frac{(\overline{n}-1)!2^{-\overline{n}}}{(\frac{\overline{n}}{2}-2)!}
  +\frac{\overline{n}!2^{-(\overline{n}+1)}}{(\frac{\overline{n}}{2}-1)!}\Big). \nonumber\\
\end{eqnarray}
Recall the Einstein-Hilbert action for manifolds with boundary  (see \cite{Wa3} or \cite{Wa4}),
\begin{equation}
I_{\rm Gr}=\frac{1}{16\pi}\int_Ms{\rm dvol}_M+2\int_{\partial M}K{\rm dvol}_{\partial_M}:=I_{\rm {Gr,i}}+I_{\rm {Gr,b}},
\end{equation}
  where
  \begin{equation}
K=\sum_{1\leq i,j\leq {n-1}}K_{i,j}g_{\partial M}^{i,j};~~K_{i,j}=-\Gamma^n_{i,j},
\end{equation}
and $K_{i,j}$ is the second fundamental form, or extrinsic
curvature. Take the metric of (3.44), then by Lemma A.2 in \cite{Wa3},
$K_{i,j}(x_0)=-\Gamma^n_{i,j}(x_0)$ when $i=j<\overline{n}+2$, otherwise is zero. Then
  \begin{equation}
K(x_0)=\sum_{i,j}K_{i.j}(x_0)g_{\partial M}^{i,j}(x_0)=\sum_{i=1}^{\overline{n}+2}K_{i,i}(x_0)=-\frac{\overline{n}+1}{2}h'(0).
\end{equation}
Hence  we have proved get the equivariant Kastler-Kalau-Walze type theorems with torsion on compact manifolds   with  boundary.
\begin{thm}
 Let M be a $\overline{n}+2$-dimensional compact spin manifold  with the boundary $\partial M$, for equivariant twisted Bismut Laplacian with torsion,
\begin{align}\label{eq:3.19}
 &\widetilde{{\rm Wres}}\big[\pi^+H_{X}^{-1}\circ \pi^+H_{X}^{\frac{(-\overline{n}+2)}{2}}\big]\nonumber\\
=&\frac{(n-2)\pi^{\frac{n}{2}}}{(\frac{n}{2}-1)!}\int_{M}\mathrm{dim}(F)2^{n}
 \Big(-\frac{1}{12}s+\frac{1}{2}\mathrm{div}^{g}(X)+g(Y,X) {\rm{Tr}}^{F}(\Phi)\nonumber\\
&+\sum _{1\leq\alpha<\beta<\gamma\leq n} 2 T^{2}_{\alpha \beta \gamma}  {\rm{Tr}}^{F}(\Phi^{2}) +\frac{1}{2}|Y|^{2}
  {\rm{Tr}}^{F}(\Phi^{2}) \Big)  \mathrm{dvol}_{M}\nonumber\\
 &+\frac{(3\overline{n}-6)2^{(\frac{\overline{n}}{2}-2)}\pi
  Vol(S_{\overline{n}})}{(\frac{\overline{n}}{2}+1)!}
  \Big(\frac{(\overline{n}-1)!2^{-\overline{n}}}{(\frac{\overline{n}}{2}-2)!}
  +\frac{\overline{n}!2^{-(\overline{n}+1)}}{(\frac{\overline{n}}{2}-1)!}
  \Big)\int_{\partial_{M}}(X_{n}-2 Y_{n} ){\rm dvol}_{\partial_{M}}
\nonumber\\
 & -\frac{(\overline{n}-2)\pi i}{(\overline{n}+1)(\frac{\overline{n}}{2}+2)!}
 \frac{(-\frac{1}{2}i\overline{n}^{4}-i \overline{n}^{3}+\frac{5}{4}i\overline{n}^{2}+\frac{5}{2}i\overline{n} )
(\frac{\overline{n}}{2}+2 )(\frac{\overline{n}}{2}+3 )\cdots (\overline{n}-1)
}{2^{\frac{\overline{n}}{2}+1}}Vol(S_{\overline{n}})\int_{\partial_{M}}K{\rm dvol}_{\partial_{M}},
\end{align}
 where $Vol(S_{\overline{n}})$ is the canonical volume of $S_{\overline{n}}$ and $s$ denotes the scalar curvature.
 \end{thm}

\section*{ Acknowledgements}
The first author was supported by NSFC. 11501414. The second author was supported by NSFC. 11771070.
 The authors also thank the referee for his (or her) careful reading and helpful comments.

\end{document}